\documentclass[11 pt]{amsart}
\usepackage{color}
\usepackage{amssymb,amsmath,amsthm,latexsym,mathrsfs}

\usepackage{hyperref}

\usepackage{times}
\usepackage[latin1]{inputenc}
\usepackage[T1]{fontenc}
\usepackage[dvips]{graphicx}

\def\N{{\mathbb N}}
\def\R{{\mathbb R}}
\def\g{{\gamma}}

\def\s{{\sigma}}

\def\d{{\delta}}

\def\p{{\prime}}

\def\B{ {\mathcal{B} } }
\def\A{ {\mathcal{A} } }

\def\F{ {\mathcal{F} } }

\def\+R{+_{_{ \!\! \R}}}

\def\phi{{\varphi}}
\def\ve{{\varepsilon}}
\def\e{{\epsilon}}

\def\what{\hat}
\def\wt{\widetilde}

\def\mR{{\mathscr{R}}}
\def\mS{{\mathscr{S}}}
\def\mT{{\mathscr{T}}}

\def\grad{\Lambda}

\DeclareMathAlphabet{\mathpzc}{OT1}{pzc}{m}{it}
\def\nn{\nonumber}
\def\nl{\vglue0.3truemm\noindent}
\numberwithin{equation}{section}

\usepackage{mathrsfs}

\DeclareMathSymbol{\Gamma}{\mathalpha}{letters}{"00}
\DeclareMathSymbol{\Theta}{\mathalpha}{letters}{"02}
\DeclareMathSymbol{\Lambda}{\mathalpha}{letters}{"03}
\DeclareMathSymbol{\Omega}{\mathalpha}{letters}{"0A}
\definecolor{cr}{rgb}{.8,0,0}

\begin{document}

\newtheorem{theo}{Theorem}[section]
\newtheorem*{theo*}{Theorem}

\newtheorem{pro}[theo]{Proposition}
\newtheorem{lem}[theo]{Lemma}
\newtheorem{cor}[theo]{Corollary}

\theoremstyle{remark}
\newtheorem{rem}[theo]{Remark}

\theoremstyle{definition}
\newtheorem{defin}[theo]{Definition}

\setcounter{tocdepth}{1}
\setcounter{secnumdepth}{1}

\title{Space-Time resonances and the null condition 
\\
for (first order) systems of wave equations
}

\author{Fabio Pusateri
\and
Jalal Shatah
}


\begin{abstract}
In this manuscript we prove global existence and linear asymptotic behavior of small solutions to nonlinear wave equations.
We assume that the quadratic part of the nonlinearity satisfies a non-resonance condition which is a 
generalization of the null condition given by Klainerman \cite{K1}.
\end{abstract}

\maketitle

\markboth{FABIO PUSATERI AND JALAL SHATAH}{SPACE-TIME RESONANCES AND THE NULL CONDITION} 


\section{Introduction}
Global existence and asymptotic behavior of small solutions to nonlinear wave equations 
has been a subject under active investigation for over fifty years.  One area of research, where 
much progress has been made, focuses on  identifying 
nonlinearities that  lead to global solutions for small initial data.  
In this manuscript we consider  first order systems  on $\R\times\R^3$ , of the form  
 \begin{equation}\tag{W}\label{w}
\left\{ 
\begin{array}{l}
\partial_t u  =i \grad u   
		+  Q_1 (u,v)   +  R_1(u,v)
\\
\partial_t v= - i \grad v +Q_2 (u, v)  +   R_2(u,v)
\\
u (1,x) = u^0 (x) \, ,  v (1,x) = v^0 (x),
\end{array}
\right.
\end{equation}
where $\Lambda := |\nabla|$, $Q_i(u,v)$ are bilinear  in $(u,v)$ and their complex conjugates, 
and $R_i$ are of degree $3$ or higher.    

In this paper we focus on determining some general conditions, naturally arising from the  {\it space time resonance analysis}, 
that guarantee global existence and scattering.  Our non-resonant condition imposed on the $Q_i$,  
roughly states that time resonant wave interactions should be limited to waves with different group velocities (spatially non-resonant waves).

Since  cubic and higher order terms   do not require any condition to ensure global existence, 
we will drop the $R_i$'s from any further consideration.  Moreover by introducing the notation  for   bilinear pseudo-product operator
\begin{align*}
T_{m(\xi,\eta)} (f,g) & := \F^{-1} \int  m(\xi,\eta) \what{f}(\eta) \what{g}(\xi - \eta) \, d\eta ,
\end{align*}
where $\hat g = \mathcal{F}g$ is the   Fourier transform of $g$,  and without any  loss of generality,   
we reduce  the system  to a single scalar equation 
\begin{equation}
\label{scalareq}
\left\{ 
\begin{array}{l}\displaystyle
\partial_t u - i \grad u = T_{q_{+,+} (\xi,\eta)} (u, u)
		+  T_{q_{-,+} (\xi,\eta)} (\bar{u}, u) + T_{q_{-,-} (\xi,\eta)} (\bar{u}, \bar{u})
\\
\\
u (1,x) = u^0 (x) \, .
\end{array}
\right.
\end{equation}
with quadratic nonlinearities.  Here and throughout the paper $+,-$ stands for $u$ and $\bar u$ respectively.

To motivate our work, we start by recalling recent applications of the space time resonance method to several problems.  
This method was introduced in \cite{GMS1,GMS3} where non resonant nonlinearites  were treated  for Schr\"odinger  equations, 
which corresponds to \eqref{w} with $\Lambda = |\nabla|^2= -\Delta$. 
In these works,    most of  existing results on  global existence and scattering    
of  small solutions were reproduced and explained by studying space time resonant frequencies.  
Subsequently  the method was applied to gravity water waves \cite{GMS2}, which corresponds to $\Lambda = |\nabla|^{1/2}$, 
and  to capillary waves \cite{GMS4}, which corresponds to $\Lambda = |\nabla|^{3/2}$.   
Thus it is natural to us to apply this method to  system \eqref{w}, where $\Lambda = |\nabla|$,  
which can be reduced to a system of nonlinear wave equations.   Our main result is:
  \newtheorem*{ntha}{Theorem A}
\begin{theo}
\label{maintheo}
Assume that system  \eqref{scalareq} is non-resonant  in the sense of definition \ref{null},  and that
the initial datum satisfies\footnote{
See the remark at the end of section \ref{secresonances} for some comments about these initial conditions.
}
\begin{align}
\label{initdatum}
{\| x u_0 \|}_{H^2} + {\left\| \grad x^2 u_0 \right\|}_{H^1} + {\left\| u_0 \right\|}_{H^N} \leq \e 
\end{align}
for some large enough integer $N$.
Then, if $\e$ is small enough, there exists a unique global solution to \eqref{scalareq} with
\begin{equation*}
{\| u(t) \|}_{L^\infty} \lesssim \frac{\e}{t} \, .
\end{equation*}
Moreover, $u(t)$ scatters  in $H^2$ to a linear solution as $t \rightarrow \infty$.
\end{theo}

Our non resonant condition defined in  \ref{null}  
turns out to include the classical null condition for wave equations \cite{K1},  
wave equations which are not invariant under the full Lorentz group,
as well as other systems  where global existence and asymptotic behavior of small solutions was not known.

\subsection{Background}  Since our system can be reduced to nonlinear wave equations, we give a 
 brief review of some of the main results about the long time existence
of solutions for systems of quadratic nonlinear wave equations on $\R^{1+3}$:
\begin{equation}
\label{quadraticsystem0}
\Box u_i = \sum a_{i, \alpha \beta}^{jk} \partial^\alpha u_j \partial^\beta u_k  \,\, + \,\, \mbox{cubic terms} 
\end{equation}
where $i = 1,\dots,N$ for some $N \in \N$, and the sum runs over $j,k=1,\dots,N$,
and all multi-indices $\alpha,\beta \in \N^4$ with $|\alpha|,|\beta| \leq 2$, $|\alpha| + |\beta| \leq 3$,
with the usual convention that $\partial_0 = -\partial^0=\partial_t$.  
Let us first recall that in $3$ space dimensions general quadratic nonlinearities have long range effects:  the $L^2$ norm of the nonlinearity, computed on a linear solution, 
decays at the borderline non-integrable rate of $t^{-1}$.
Thus, quadratic nonlinearties can contribute to the long time behavior of solutions.  
It is in fact known since the pioneering works of  John \cite{John0,John2}
that  finite time blowup  can occur even for solutions with small data.
On the other hand, for some very general classes of quadratic nonlinearities
solutions
were shown to exist 
and almost globally 
by John and Klainerman \cite{JK} and Klainerman \cite{K0}.

The main breakthrough in identifying   classes  of nonlinear  wave equations  
where solutions with small data exist  globally and  scatter   was in the works of   Klainerman \cite{K1},  
Choquet-Bruhat and Christodoulou \cite{C2},  and Christodoulou \cite{C1}. 
The class of nonlinearities that satisfy the ``null condition'' 
were  introduced   by Klainerman \cite{K1}, and for semlinear systems
\begin{equation}
\label{quadraticsystem}
\Box u_i = \sum_{|\alpha|,|\beta| = 1}   
      a_{i, \alpha \beta}^{jk} \partial^\alpha u_j \partial^\beta u_k  +  \,\, \mbox{cubic terms}  
\end{equation}
are given by the condition
\begin{equation}
\label{nullcondition}
\sum
  a_{i, \alpha \beta}^{jk} \xi_\alpha \xi_\beta = 0 
  \qquad \mbox{for any $\xi \in \R^4$ such that $-\xi_0^2 + \xi_1^2 + \xi_2^2 + \xi_3^2 = 0$} \, .
\end{equation}
For such systems it was shown by Klainerman \cite{K1} that in $3+1$ dimension small data solutions exists globally.  This seminal work of  Klainerman    is based on the invariance of Minkowski space under the Lorentz group and on energy estimates using the vector fields that generate the Lorentz group \cite{K0}.



Later on, building on Klainerman's original ideas, the problem of bypassing
the use of the full invariance under the Lorentz group was dealt with by other authors.
In \cite{KS} Klainerman and Sideris 
proved almost global existence of solutions for quadratic systems
\eqref{quadraticsystem0} in divergence form, under the sole assumption of translation, rotation and scaling invariance.
Further developments were made by Sideris in \cite{sid1,sid2}, where global existence of nonlinear elastic waves 
is proven under the assumption of the null condition. 
Similar results include the almost global existence of solutions
contained in the works of Keel, Smith and Sogge \cite{KSS1,KSS2}.
It also worth mentioning that several works have dealt with the question of identifying other conditions  (weaker than the null condition) 
under which global existence of solutions of \eqref{quadraticsystem0} can be proven; 
see for instance Lindblad \cite{LindWave}, Alinhac \cite{Alinhac1}, Lindblad and Rodnianski \cite{LR1,LR2}, and Katayama \cite{Kata}.

Another  approach that identifies the effects of nonlinearities on the long time  behavior of solutions is based on time resonant computations.
For ODE's this  is the Poincare-Dulac normal form.   
For PDE's  normal forms  were used by Shatah \cite{ShatahKGE}
and Simon \cite{Simon} who treated, respectively, the Cauchy problem and the final state problem 
for the Klein-Gordon equation in $3+1$ dimensions (see also Kodama \cite{Kodama} for an early appearance of normal forms in the PDE setting).   
Similar results were obtained by Klainerman using the vector fields method \cite{KKG}.

In the past several years a new algorithmic  method, called the ``space-time resonance method'',  
was developed by Germain, Masmoudi, and Shatah,
to study long time behavior of spatially localized small solutions to dispersive equations.  
By bringing together ideas from both vector fields and normal forms,
this new method proved to be effective in proving new results \cite{GMS2,Ger,GM} 
as well as simplifying already existing ones \cite{GMS1,GMS3,KP}.  
A description of this method can be found in \cite{GMS1}.

\subsection{Notations}   
We use $R$ to denote indistinctly any one of the components of the vector of
Riesz transforms $R = \frac{\nabla}{\Lambda}$, where $\Lambda := |\nabla|$.   
 $L^p$ norms will be denoted either by ${\| \cdot \|}_{L^p}$ or ${\| \cdot \|}_{p}$.
For $s \geq 0$, $p \geq 1$, we define the usual Sobolev norms
\begin{align*}
{\|\varphi\|}_{W^{s,p}} & := {\| \langle \nabla  \rangle^{s} \varphi\|}_{L^p}
\, ,
\\
{\|\varphi\|}_{\dot{W}^{s,p}} & := { \| {\grad}^{s} \varphi \|}_{L^p}
\, ,
\end{align*}
where $\langle x \rangle := {(1+{|x|}^{2})}^{1/2}$.
We let $H^s := W^{s,2}$ and $\dot{H}^s := \dot{W}^{s,2}$.
\nl
Finally we write $A \lesssim B$ to mean $A \leq C B$ for some positive absolute constant $C$.

\section{Resonance analysis and non-resonant bilinear forms}

\subsection{Resonance analysis}
Recall that to compute resonance for an equation of the type
\begin{equation*}
i u_t   +  P(\frac 1i\nabla) u  =  T_{m(\xi,\eta)} (u,u),
\end{equation*}
we write  Duhamel's formula in Fourier space for the ``profile'' of $u$, namely  $f:=e^{-itP(\frac 1i\nabla)}  u$, 

\begin{equation}
\label{Duhamel}
\what{f}(t,\xi) = \what{u}_0 (\xi) + \int_0^t \!\! \int e^{is \phi (\xi,\eta)}
		m(\xi,\eta)  \what{f} (s,\eta) \what{f} (s,\xi - \eta) \, d\eta \, ds \, ,
\end{equation}
where $ \phi (\xi,\eta) := - P(\xi) + P(\eta) + P(\eta - \xi)$  (with obvious signs modifications occur if $Q$ depends also on $\bar{u}$), and define the {\it time resonant set} 
\begin{equation*}
\mT = \{  (\xi,\eta) \, : \,  \phi (\xi,\eta) = 0  \} \qquad \mbox{(no oscillations in $s$)} \, , 
\end{equation*}
 the {\it space resonant set} 
\begin{equation*}
\mS = \{  (\xi,\eta) \, : \, \nabla_\eta \phi (\xi,\eta) = 0  \} \qquad \mbox{(no oscillations in $\eta$)} \, , 
\end{equation*}
and the {\it space-time resonant set} 
\begin{equation*}
\mR =  \mT  \cap  \mS \, .
\end{equation*}

Since for system \eqref{scalareq}  both $u$ and $\bar u$ are present in the bilinear terms, there are three types of interactions that we need to analyze.

\subsubsection{The - - case}
The phase $\phi_{- -} = -|\xi| - |\eta| - |\xi -\eta|$ clearly vanishes only at $\xi = \eta = 0$:
\begin{equation*}
\mT_{- -} =  \{ \eta = \xi = 0 \} \, .
\end{equation*}
Since the time resonant set is reduced to a point, we can perform a {\it normal form transformation}.
This allows us to obtain the $L^\infty$ decay in a more direct fashion
(without the need to resort to weighted estimates).
For completeness we compute
\begin{equation*}
  \mathscr S_{--}=\{\eta= \lambda \xi, \quad 0\le \lambda \le 1\} ,
\end{equation*}
and
\begin{equation*}
\mR_{- -} =  \{ \xi = \eta = 0 \} \, .
\end{equation*}

\subsubsection{The + + case}
\label{secresonances}
The phase $\phi_{+ +} = -|\xi| + |\eta| + |\xi -\eta|$ 
vanishes on 
\begin{equation*}
\mathscr T_{++}= \{\eta= \lambda \xi, \quad 0\le \lambda \le 1\} \, .
\end{equation*}
A simple computation shows that
\begin{equation*}
  \mathscr S_{++}=\{\eta= \lambda \xi, \quad 0\le \lambda \le 1\} \, ,
\end{equation*}
whence
\begin{equation*}
\mR_{+ +} = \mS_{+ +} = \mT_{++} \, .
\end{equation*}
The space time resonant set is very large, thus some additional structures are  needed to help in controlling these resonances.    
The first  structure will be imposed on the interaction by requiring   the symbol $q_{++}$  to vanish on $\mR_{+ +}$.  
The second structure is present in the phase
\[
|\xi|\nabla_\xi\phi_{++} = \frac {\eta-\xi}{|\eta-\xi|}\phi_{++}   -|\eta|\nabla_\eta\phi_{++} \, ,
\]
and can be interpreted by saying that all resonant waves have the same group velocity and thus are spatially localized in the same region.   
This  fact together with  
 $\mS_{+ +} =  \mT_{+ +}$ allows to control these resonances by  solely relying  on weighted energy estimates.

\subsubsection{The - + and + - cases}
Up to the change of variables\footnote{ 
notice that $m_{12}(\xi, \xi - \eta) = - m_{12}(\xi, \eta)$.
}
$\eta \rightarrow \xi - \eta$, these two case are the same.
Therefore, we will just focus on the $- +$ case.
Since the phase is  $\phi_{- +} = -|\xi| - |\eta| + |\xi -\eta|$, then
\begin{align*}
    & \mathscr T_{-+}=  \{\eta= \lambda \xi, \quad  \lambda \le 0\}\cup \{ \xi = 0 \} \, , 
\\
    &  \mathscr S_{-+}= \{\eta= \lambda \xi, \quad \lambda \le 0 \quad \mbox{or} \quad \lambda \ge 1\}\cup \{ \xi = 0 \} \, ,  
\\
   &  \mathscr R_{-+} = \{\eta= \lambda \xi, \quad  \lambda \le 0\}\cup \{ \xi = 0 \}. 
\end{align*}

Again the set $\mR_{-+}$ is  very big  and additional conditions are needed to ensure global existence 
and linear asymptotic behavior of solutions.  
These conditions are similar to the $++$ interaction, i.e., $q_{-+}$ vanishes on $\mR_{-+}$ and the fact that

\[
|\xi|\nabla_\xi\phi_{-+} = \frac {\eta-\xi}{|\eta-\xi|}\phi_{-+}  + |\eta|\nabla_\eta\phi_{-+}.
\]

However this interaction presents an additional difficulty over the $++$ case since $\mT_{-+}\subsetneq \mS_{-+}$, 
which requires both normal forms transformation and weighted estimates.  The fact that this is an added difficulty is explained below.
\subsection{Non-resonant bilinear forms}

From  Duhamel's formula  for equation \eqref{scalareq} in Fourier space the quadratic term is expressed as

\begin{equation*}
\what{B}_{\pm\pm} (t,\xi) = \int_1^t \! \int  e^{is \phi_{\pm\pm}(\xi,\eta)}  q_{\pm\pm} (\xi, \eta)
			\what{f}_\pm (s,\eta) \what{f}_\pm (s,\xi-\eta)  \, d\eta ds
\end{equation*}
where $f_+ = f$ and $f_- = \bar{f}$ and
$$
\phi_{\e_1,\e_2} (\xi,\eta) =  - |\xi| + \e_1 |\xi-\eta| + \e_2 |\eta| \, ,
$$
for $\e_i = \pm$.   The quadratic interaction is given in terms of its symbol $q_{\pm\pm}$.   To define non-resonant bilinear forms we start by defining the class  of symbols that we will be dealing with,
\begin{defin}
\label{definop}
A symbol $m = m (\xi,\eta)$ belongs to the class $\B_s$ if\\[.5em]
\noindent $\bullet$  It is homogeneous of degree $s$;\\[.5em]
\noindent $\bullet$  It is smooth outside of $\{ \xi = 0 \} \cup \{ \eta = 0 \} \cup \{ \xi - \eta= 0 \}$.\\[.5em]
\noindent $\bullet$  For any labeling $(\xi_1,\xi_2,\xi_3)$ of the three Fourier variables $(\xi,\eta,\xi-\eta)$ the following holds:
\begin{align*}
\mbox{for} \quad |\xi_1| \ll |\xi_2|,|\xi_3| \sim 1 \qquad m = \A \left( |\xi_1|, \frac{\xi_1}{|\xi_1|}, \xi_2 \right)
\end{align*}
for some smooth function $\A$.
\end{defin}

\noindent
Loosely speaking, symbols in $\B_0$ are Coifman-Meyer \cite{CM}
except, possibly, along the axes $\xi=0,\eta=0$ and $\xi-\eta=0$, 
where they can have singularities like linear Mihlin-H\"{o}rmander multipliers.
Symbols in $\B_s$ are essentially of the form $|\xi|^s m_0$ for some $m_0 \in \B_0$.  The  boundedness of these bilinear operators  on $L^p$  is given in the appendix \ref{prooftheoop}.
With the class $\B_s$ defined above we can define a non-resonant  system as
\begin{defin}[\bf{Non-resonant bilinear forms}]
\label{null}
System \eqref{scalareq} is called non-resonant  if
\begin{equation}
\label{SYM}
q_{\pm\pm} (\xi,\eta) = a(\xi,\eta) \phi_{\pm\pm} (\xi, \eta) +  b(\xi,\eta) \cdot \nabla_\eta \phi_{\pm\pm} (\xi, \eta) \, ,
\end{equation}
with $a \in \B_{-1}$  and $b \in \B_0$. Additionally we require that
\begin{align}
\label{restriction1}
{|\eta|} \nabla_{(\xi,\eta)} a (\xi,\eta) \quad \mbox{or} \quad {|\xi-\eta|} \nabla_{(\xi,\eta)} a (\xi,\eta) 
	= \frac{\mu_0^{(1)}(\xi,\eta)}{|\xi|}
		+ \frac{\mu_0^{(2)}(\xi,\eta)}{|\eta|} + \frac{\mu_0^{(3)}(\xi,\eta)}{|\xi-\eta|}
\end{align}
for some $\mu_0^{(i)} \in \B_0$, 
and
\begin{align}
\label{restriction2}
{\left\| T_{b(\xi,\eta)} (f,g)  \right\|}_{L^2} \lesssim  
	    {\left\| f \right\|}_{L^2}  \sum_{j=0}^k {\left\| R^j g \right\|}_{L^\infty}
	+   \sum_{j=0}^k {\left\| R^j f \right\|}_{L^\infty} {\left\| g \right\|}_{L^2}
\end{align}
for some $k \in \N$. 
\end{defin}

Some comments about this definition.  Equation \eqref{SYM}  asserts that bilinear interaction vanishs on $\mathscr R_{\pm\pm}$, the space time resonant set.  
The presence of  $\varphi_{\pm\pm}$ in equation \eqref{SYM}  allows us to perform normal form  transformation on one part of the bilinear terms  
(integration by parts in  $s$), while the presence of  $\nabla_\eta\varphi_{\pm\pm}$ allows us to treat the remaining part by weighted estimates 
(integration by parts in  $\eta$).  The {\it classical null condition} is equivalent to $a=0$ (see below).  
Equation \eqref{restriction1} essentially avoids having  $a(\xi,\eta) \sim {|\xi|}^{-1}$ which would be too singular to handle.   
Equation \eqref{restriction2}  is needed due to  failure of
the $L^2 \times L^\infty \rightarrow L^2$ estimate for symbols in $\B_0$.
It would be possible to avoid this last technical restriction by resorting to the use of Besov spaces,
but for the sake of simplicity we do not pursue this matter here
(see also remark \ref{remaboutR} for more comments about this aspect).

\subsection{Examples of non-resonant systems}  Now we give examples of non-resonant systems and explain how they relate to  existing definitions on `` null systems", and how our definition is a natural extension of previous ones.

\subsubsection{Classical null forms.}
Quadratic (semilinear) nonlinearities satisfying the null condition \eqref{nullcondition}
are linear combinations of
\begin{subequations}
\label{null*}
\begin{align}
\label{Q_ij}
Q_{ij}(u,v) & = \partial_i u \partial_j v - \partial_j u \partial_i v \, ,
\\
\label{Q_i0} 
Q_{0i}(u,v) & = \partial_t u \partial_i v - \partial_i u \partial_t v \, ,
\\
\label{Q_0}
Q_{0}(u,v) & = \partial_t u \partial_t v - \nabla u \! \cdot \! \nabla v \, .
\end{align}
\end{subequations}
By letting $(u_\pm, v_\pm) = (\partial_t u \mp i \grad u, \partial_t v \mp i \grad v)$,
one can reduce systems \eqref{quadraticsystem} to first order systems
in the unknowns $u_\pm$ and  $v_\pm$.
Then, one can check that the respective symbols of the above null forms,
corresponding to interactions of $u_{\e_1}$ and $v_{\e_2}$, are given (up to a constant factor) by
\begin{align*}
m_{ij}^{\e_1,\e_2}  (\xi,\eta) & = 2 \frac{\eta_i \xi_j - \eta_j \xi_i}{|\eta||\xi-\eta|}
	= \partial_{\eta_i} \phi_{+ +} (\xi,\eta) \partial_{\eta_j} \phi_{+ -} (\xi,\eta) 
	-  \partial_{\eta_j} \phi_{+ +} (\xi,\eta) \partial_{\eta_i} \phi_{+ -} (\xi,\eta)
\\ 
m_{0i}^{\e_1,\e_2} (\xi,\eta) & = 
	\left( \e_1 \frac{\eta_i}{|\eta|} + \e_2 \frac{\eta_i - \xi_i}{|\eta-\xi|} \right)
	= \partial_{\eta_i} \phi_{\e_1 \e_2} (\xi,\eta)
\\
m_{0}^{\e_1,\e_2} (\xi,\eta) & =
	2 \left(  1 - \e_1 \e_2 \frac{\eta}{|\eta|} \cdot \frac{\xi - \eta}{|\xi - \eta|} \right)
	=
	{\left| \nabla_\eta \phi_{\e_1 \e_2} (\xi,\eta) \right|}^2 \, .
\end{align*}

These symbols are of the form \eqref{SYM} with $a=0$ and vanish on the space resonant set.
Thus  classical null forms are spatially non-resonant, and therefore can be treated by weighted estimates 
and without normal forms transformation.  
Note that  in our system the $+-$ interactions have $\mT_{+-}\subsetneq \mS_{+-}$,  
and our symbols \eqref{SYM} only vanish on the smaller set $\mT_{+-}$.
To treat these interactions a normal form transformation is needed, leading to  terms which are not well spatially localized.  
This causes a difficulty that will be elaborated on below.

The classical quasilinear null forms are also non-resonant 
in that their symbols satisfy \eqref{SYM} where the homogeneities of $a$ and $b$ are increased by $1$.
However, for general first order systems of the type \eqref{scalareq}, 
quasilinear equations lose derivatives in the energy estimates unless there are cancellations  present.  
In the special case where the first order quasilinear system  comes from a second order system of 
 wave equations with quasilinear null terms, cancellations are present in the energy estimates.
Thus, our results will apply for such nonlinearities as well\footnote{
More specifically, in the quasilinear case the most efficient proof of the analogue of Theorem \ref{maintheo}
would consist of two steps:  1) establishing energy and weighted energy estimates directly on the second order wave equation
so to obtain the weighted bounds in \eqref{Xnorm};
2) run our proof to show the decay of solutions.}.

\subsubsection{Systems with multiple speeds}
For systems
\begin{equation}
 \partial_t u_\ell - ic_\ell  \grad u_\ell  = \sum_{j,k}\,T_{ q^\ell_{j,k} (\xi,\eta)} (u_j, u_k)
\end{equation}
the phases are given by $-c_\ell|\xi| + c_j|\xi - \eta| +c_k|\eta|$.
In the case $c_j \neq \pm c_k$ one has $\mS= \{0,0\}$, so that our results will trivially apply.
If $c_j = \pm c_k$, $\mR = \{0,0\}$ unless $c_k = \pm c_\ell$. 
Therefore, in the case $c_k \neq \pm c_\ell$, global existence can be obtained provided a suitable null condition is imposed at $\{0,0\}$. 
This is similar to the work on quadratic NLS \cite{GMS3}.
The full non-resonance condition is then needed only for interactions of the form $ -c_\ell|\xi| + c_\ell|\xi - \eta| \pm c_\ell|\eta|$.
This extension is similar to the result of Sideris and Tu \cite{ST}. 
We also refer the reader to the work of Katayama and Yokoyama \cite{KataYoko} and references therein
for more on systems with multiple speeds. Some examples of interest that can be treated using our techniques are:

\noindent 
$1.$ First order systems of the form
\begin{equation}
\label{ex1}
\left\{ 
\begin{array}{l}\displaystyle
\partial_t u + i c \grad u = {|v|}^2
\\
\partial_t v + i \grad v = T_{m} (v,v) + u^3
\end{array}
\right.
\end{equation}
where $c>1$ and $m$ is non-resonant as in definition \ref{definop}.
Here no special null condition at the origin is assumed on the bilinear form in the first equation.    This system does not satisfy any existing null condition criteria set by the vector fields method.  In fact we believe that this   system is not amenable to analysis by the vector fields method due to the simultaneous failure of the Lorentz invariance and the need of a normal form transform.

Our method works by first applying a normal form transformation on $u$ (notice that the phase is bounded below by $(c-1)|\xi|$),  
and then handling the singularity introduced by such a transformation through a  spread-tight splitting  explained in section \ref{secproof}.  

\noindent 
$2.$ Systems of wave equations of the form
\begin{equation}
\label{ex2}
\left\{ 
\begin{array}{l}\displaystyle
\Box_1 u = {(\partial v)}^2  +  \partial u {(\partial w)}^2 
\\
\Box_1 v = \partial v \partial w + Q(v,v)
\\
\Box_c w = \partial v \partial w  +  \partial u {(\partial v)}^2
\end{array}
\right.
\end{equation}
where $\Box_a := \partial_t^2 - c^2 \Delta$, $c \neq 1$ and $Q$ is any null form.
This is an example of a nonrelativistic system satisfying the weak null condition.
Global existence can be obtained with a weaker decay on $u$ of the form
${\| \partial u \|}_{L^\infty} \lesssim t^{-1+\e}$, for $\e \ll 1$.


\subsubsection{Non-locality and absence of Lorentz invariance}\label{remnonlocal}
The class of systems \eqref{w}, under the non-resonance condition given by definition \ref{null},
includes the class of second order wave equations 
\begin{equation}
\label{nonlocalnull}
\Box u = T_1 Q (T_2 u, T_3 u)
\end{equation}
where the  $T_i$'s are zero-th order operator and $Q$ is any combination of the nonlinear terms \eqref{Q_ij}-\eqref{Q_0}
(or their quasilinear version).
For systems as \eqref{nonlocalnull} the action of the Lorentz boosts $L_i = x_i \partial_t + t \partial_i$ on the nonlinearity
produces some terms which are too singular to be estimated.
This makes the classical \cite{K1} vector fields method difficult to apply.

In \cite{sid1,sid2} Sideris considered quasilinear hyperbolic systems governing
the motion of isotropic, homogeneous, nonlinear elastic waves. 
Like systems with multiple speeds, these systems are only classical invariant, i.e. they do not possess Lorentz invariance.
By imposing a null condition on the nonlinear terms of the form $F(\nabla u) \nabla^2 u$,
he was able to show global existence of solutions. 
As mentioned in the introduction, several other works have dealt with the problem of long time existence
for classically invariant systems on $\R^{3+1}$, see for example \cite{KS,KSS1,KSS2}.
Our methods are also applicable to the systems considered in these works.


\subsubsection{A remark about the initial data}\label{remdata}
In contrast with the results mentioned previously, 
our initial data belong to a low weighted Sobolev space. 
In particular we only ask for $x u_0 \in H^2$ and ${|x|}^2 \grad u_0 \in H^1$, see \eqref{initdatum}.
In comparison, the spaces used in \cite{KS,sid1,KSS1} would require
${( |x| \grad )}^i u_0 \in \grad (L^2)$, for $i= 0 \dots k$ and some $k \geq 7$.
This means that we can allow more oscillating data.
For example, for data behaving at infinity like $\cos |x| / {|x|}^\alpha$, we can allow any $ \alpha > \frac{7}{2}$, 
whereas in the other works one would need $\alpha > \frac{17}{2}$.



\section{Outline of the proof}\label{secproof}

Before we outline the  proof of Theorem \ref{maintheo} we would like to point out two difficulties in our problem: 
\nl
$a)$ Although  the space-time resonance method is algorithmic, 
its implementation  on the space-time resonant set is very much problem dependent.  
This is  due to the fact that the aforementioned set can be large with no clear  criteria, set by the method, 
to address  how large  is large.    
Its application to nonlinear dispersive equations  has been restricted so far to cases where the resonant set is very small.
In particular, for problems such as the Schr\"odinger equation, the resonant set is a point; 
and for gravity water waves, there are no quadratic  time resonances.
However for hyperbolic systems this set is large.    
For the system we are considering here  the space time resonant set is  
$4$ dimensional in a $6$ dimensional space.   
Treating such a big space-time resonant set required new ideas, which we present in this manuscript.
\nl
$b)$ When  space resonant frequencies are different from time resonant frequencies, and when both types of resonances are present 
in the bilinear interactions,  a normal form transformation is needed
\begin{equation*}
u \to  u + T_m(u,u)   \Leftrightarrow \what{f} \to  
\what{f} +  \int e^{is \phi (\xi,\eta)} m(\xi,\eta)  \what{f} (s,\eta) \what{f} (s,\xi - \eta) \, d\eta \, .
\end{equation*}
The bilinear term $T_m(u,u) $  need not be well localized in space since 
the outcome of the interaction may have a different group velocity, i.e., $\nabla_\xi \phi \ne 0$, which is the case  
for  \eqref{scalareq}. 
Thus weighted estimates on this bilinear interaction tend to grow at a fast rate with time.  
We refer to these bilinear interactions  as  {\it spread terms}.  
This is in contrast to non space resonant frequencies which are well spatially localized and 
for which weighted estimates tend to grow very slowly. 
We refer to such bilinear interactions as {\it tight terms}. The presence of tight and spread terms  
is  problematic and requires a careful analysis when trying to establish the decay of solutions.
This is the case here for the $+-$ interactions, as was the case for the 2D Schr\"odinger equation \cite{GMS3}.
Our strategy in obtaining the pointwise decay of solutions will be explained below in more details.



\subsection{Reduction of system $\boldsymbol{\eqref{scalareq}}$} 
By isolating the terms in equation \eqref{scalareq} that are most difficult to estimate, 
we can considerably  simplify our notation and presentation.

\subsubsection{ Reduction to the $- +$ case}
As the analysis of resonances in section \ref{secresonances} showed,
the  $- +$ interactions lead to a more complicated resonant set than  the  $+ +$ and $- -$ interactions.
The $- +$ case actually contains the difficult aspects of both the $+ +$ and $- -$ cases.
More precisely, in the $-+$ case we will need  to decompose the phase space in two sets:
one containing $\mR_{- +}$ but not $\mS_{- +} \cap \mT_{- +}^{c}$,
and the other one containing $\mS_{- +} \cap \mT_{- +}^{c}$ but not $\mR_{- +}$.
The analysis on the set containing $\mR_{- +}$ and  $\mS_{- +} \cap \mT_{- +}^{c}$,
would be enough to take care of the $+ +$ and the $- -$ interaction, respectively.


Therefore, from now on we focus only on this type of interaction, 
and we will drop the $- +$ indices for lighter notations.

\subsubsection{Reduction to $a(\xi,\eta)=\frac{1}{|\eta|}$ and $b(\xi,\eta)=1$}
Recall that we are imposing the restriction \eqref{restriction1} on $a$. 
This means that $a$ can have singularities of the type $1/|\eta|$ or $1/|\xi-\eta|$, but not of the form $1/|\xi|$.
By the symmetry between $\eta$ and $\xi-\eta$, 
we can then assume that $a$ is of the form ${\mu_0 (\xi,\eta)}/{|\eta|}$ for some $\mu_0 \in \B_0$.
Moreover, since the presence of symbols in the class $\B_0$ is irrelevant for our estimates
on the terms corresponding to the symbol $a(\xi,\eta) \phi(\xi,\eta)$,
we can simply assume that $a$ is given by $1/|\eta|$.

Finally, since we assume that $b$ satisfies  \eqref{restriction2}, 
and since we will show that $\|R^j u\|_{L^\infty}$ is controlled with a decay of 
$t^{-1}$ (see remark \ref{remaboutR} below) 
we can reduce matters to $b = 1$.
It will be clear to the reader what minor modifications are necessary to 
perform the estimates for a general $b \in \B_0$ and satisfying \eqref{restriction2}.

In view of these reductions the non-resonant equation becomes
\begin{equation}
\label{Duhamelf}
\what{f} (t,\xi) =  \what{f}_0 (\xi) +
\int_1^t  \int e^{is\phi(\xi,\eta)} \left( 
	\frac{\phi(\xi,\eta)}{|\eta|} 
	 + \nabla_\eta \phi(\xi,\eta)\right)
	\what{f}(s,\xi-\eta) \what{f}(s,\eta) d\eta ds \, .
\end{equation}
Furthermore, we recall that $\nabla_\xi \phi$ vanishes on the resonant set, and in particular the following identity holds:
\begin{equation}
\label{xidxiphi}
| \xi | \nabla_\xi \phi = \frac{\eta-\xi}{|\eta-\xi|} \phi + |\eta| \nabla_\eta \phi \, .
\end{equation}

\subsection{Splitting of the profile  $f$}
Integrating by parts  in $s$  in the terms containing the phase $\phi$ we  get 
\begin{equation*}
\what f (t,\xi)  \overset{def}{=} \what f_0 (\xi) +  \what g (t,\xi) +  \what h (t,\xi) \stackrel{def}{=} \what f_0 (\xi) +  \what g (t,\xi)  +  \what{h}_0 (t,\xi) + \what{h}_1 (t,\xi)
\end{equation*}
where
\begin{subequations}\label{fgh}
\begin{align}
\label{ghat}
\what{g} (t,\xi)& \stackrel{def}{=} \int e^{it\phi(\xi,\eta)}  
			\frac{1}{|\eta|} \what{f}(t,\xi-\eta) \what{f}(t,\eta) \, d\eta \, ,\\
\label{h0hat}
\what{h}_0 (t,\xi) & \stackrel{def}{=} \int_1^t \! \int e^{is\phi(\xi,\eta)} 
	\nabla_\eta \phi(\xi,\eta)  \what{f}(s,\xi-\eta) \what{f}(s,\eta) \, d\eta ds,
\\
\label{h1hat}
\what{h}_1 (t,\xi) & \stackrel{def}{=} \int_1^t \! \int e^{is\phi(\xi,\eta)} 
		\frac{1}{|\eta|} \partial_s \left( \what{f}(s,\xi-\eta) \what{f}(s,\eta) \right) \, d\eta ds \, .
\end{align}
\end{subequations}

This splitting can be understood in the following manner:  1) $g$ comes from the normal form transformation, has very good time decay but  is spatially  spread; 2) $h_0$ is a spatially tight term due to the presence of $\nabla_\eta\phi(\xi,\eta)$;  and 3)   $h_1$ is a cubic in  $f$.

Regarding  $h_0 = h_0 (f,f)$ as  a bilinear form of $f= f_0 + g + h$,  
we can decompose $h_0$ 
\begin{align*}
h_0 (f,f) &= h_0 (f_0, f_0 + h) + h_0 (f_0, g) + h_0 (g, f) + h_0 (h, h)  +  h_0 (h, g) 
		 +  h_0 (h, f_0) \\
		&= h_0 (g, f) + h_0 (h, h)  + h_*,
\end{align*}
and thus decompose $f$ further
\begin{equation}
\label{decompf}
f= f_0 + g + h_0(h,h) + h_0(f,g) + h_1 + h_*.  
\end{equation}

\subsection{Strategy of the proof and organization of the paper}
The proof of global existence will follow from the following a priori bounds on $u = e^{it\grad} f$ :
\begin{equation}\label{apriori}
\left\{ 
\begin{array}{lll}
 { \| u \| }_{H^N} \lesssim t^\ve, &  { \| u \| }_{H^2} \lesssim 1, &{ \| u \| }_{L^\infty}, {\| R u \|}_{L^\infty} 
		\lesssim \frac{1}{t},\\
 { \| x f \| }_{L^2} \lesssim t^\g, &  { \| \grad x f \| }_{H^1} \lesssim 1,  
	& { \left\| {|x|}^2 \Lambda f \right\| }_{H^1} \lesssim t  \, ,
\end{array}
\right.
\end{equation}
and a continuation of the local-in-time solution. Here, $N$ is a suitably large number  and $\ve$ and $\g$ are arbitrarily small fixed  positive constants.  This leads us to 
define  the space $X$ by the norm associated to these bounds:
\begin{equation}\label{Xnorm}
\begin{split}
\| u \|_X :=  \sup_{t\geq 1}& \left[   t^{-\ve} { \| u \| }_{H^N} + \| u \|_{H^2}
+t (   {\| u \|} _{L^\infty}  + {\| R u \|}_{L^\infty}  ) \right.  \\
& \left. +  t^{-\g} { \| x f \| }_{L^2} + { \| \grad x f \| }_{H^1} + t^{-1} { \left\| {|x|}^2 \Lambda f \right\| }_{H^1}
\right] \, .
\end{split}
\end{equation}

\begin{rem}
\label{remaboutR}
The presence of ${\| Ru \|}_{L^\infty}$ is 
not surprising because the Riesz transform
is already present in the interaction symbol $\nabla_\eta \phi$.
However we remark here that if $b$ is any symbol in $\B_0$ satisfying \eqref{restriction2},
the same $X$-norm above would work.
Indeed, as a byproduct of our estimates we have
\begin{equation}
\label{RkuLinfty}
{\left\|  R^k u \right\|}_{L^\infty} \lesssim \frac{1}{t} \left[ \e + {\| u \|}^2_X \right] 
\end{equation}
for any $k$.
This is because $L^\infty$ estimates on $u = e^{-it\grad} u_0  + e^{-it\grad} g + e^{-it\grad} h$
are obtained by
\begin{itemize}
\item[a)] using Sobolev's embedding on $g$: 
${ \| e^{-it\grad} g \| }_{L^\infty} \lesssim { \| e^{-it\grad} g \| }_{W^{1,p}}$ for $p \gg 1$,
and then showing $t { \| e^{-it\grad} g \| }_{W^{1,p}} \lesssim {\| u \|}_X^2$;
\item[b)] estimating weighted $L^2$ norms of the main components of $h$ by means of \eqref{disp0} or \eqref{disp},
and using the same argument in $a)$ on the remaining components.
\end{itemize}
In both of these operations the presence of Riesz transforms becomes irrelevant.
\end{rem}

A key aspect of our proof is the different treatment of the components $g$ and $h$,
and the different treatment of some components of $h$ itself.
In particular, the bound on ${\|u\|}_X$ will follow from the following bounds on $g$ and $h$:
\begin{align}\label{boundsg}&\left\{
\begin{array}{l}
 { \| g \| }_{H^N} \lesssim {\| u \|}_X^2, \quad 
{ \| x g \| }_{L^2} \lesssim t^\g {\| u \|}_X^2,  \quad
{ \| \grad x g \| }_{H^1} \lesssim {\| u \|}_X^2 ,
\\[.5em]
{ \left\| \grad {|x|}^2 g \right\| }_{H^1} \lesssim t {\| u \|}_X^2, \quad 
{ \left\| e^{it\grad} g \right\| }_{L^\infty} \lesssim\frac 1t {\| u \|}_X^2 \, .
\end{array}
\right.\\[1em]
\label{boundsh}&\left\{
\begin{array}{l}
{ \| h \| }_{H^N} \lesssim  t^\ve  {\| u \|}_X^2, \quad
{ \| h \| }_{H^2} \lesssim {\| u \|}_X^2, \quad
{ \| x h \| }_{L^2} \lesssim t^\g {\| u \|}_X^2, \quad
{ \| \grad x h \| }_{H^1} \lesssim {\| u \|}_X^2, 
\\[.5em]
{ \left\| \grad {|x|}^2 h \right\| }_{L^2} \lesssim t^a {\| u \|}_X^2, \quad
{ \left\| \grad^2 {|x|}^2 h \right\| }_{L^2} \lesssim t^b {\| u \|}_X^2, \quad
{ \left\| e^{it\grad} h \right\| }_{L^\infty} \lesssim \frac{1}{t} {\| u \|}_X^2.
\end{array}\right.
\end{align}
where $a$ and $b$ are (small) positive constants  satisfying
$0 < \g < b < \frac{a}{3}$ and  $ a < \frac{1}{8}$.

These a priori bounds will imply global existence provided the data is small:
\begin{equation}
\label{quadraticestimate}
{\left\| e^{it\Lambda}(g+h) \right\|}_X \lesssim {\| u \|}^2_X  \implies {\| u \| }_X \lesssim  {\left\| e^{it\Lambda}f_0\right\|}_X   +  {\left\| e^{it\Lambda}(  g+h )\right\|}_X 
\lesssim \e + {\| u \|}^2_X \, ,
\end{equation}
which in turn gives ${\| u \|}_X \lesssim \e$.

From \eqref{apriori} and \eqref{boundsh} we see that $h$ has the same energy and pointwise estimates as $f$,
and better weighted estimates than $f$.  
Thus, the bilinear terms that need to be bounded are $g$, $h_0 (h, h)$, $h_0 (f, g)$ and  $h_1$.   
All the remaining bilinear terms, denoted by  $h_*$ in \eqref{decompf},  
are easier to estimate because their arguments satisfy stronger bounds.

In what follows we briefly describe the organization of the paper together with the main steps needed in the proof.
Estimates for the Sobolev norms are pretty straightforward, since we are dealing with a semilinear equation. 
These are shown in section \ref{secenergy}.

In section \ref{weightedg} we obtain weighted and $L^\infty$ estimates for the spread component $g$.
This component is the one whose responsible for the fast growth in time of the weighted norms.
On the other hand, since it consists of a bilinear term with no time integration, 
its decay in $L^\infty$ can be obtained very easily.

In section \ref{secweightedh_0} we prove a priori bounds on weighted $L^2$ norms of $h_0(h,h)$. 
Thanks to the presence of the symbol $\nabla_\eta \phi$, and to the identity \eqref{xidxiphi}, 
we can always integrate by parts at least twice in time and/or frequency .
As a consequence we can prove that $h_0(h,h)$ satisfies the stronger weighted bounds \eqref{boundsh} 
that hold for the  $h$ component.

Section \ref{Linftyh} contains the $L^\infty$ estimate for $h_0 (h,h)$.
In order to obtain the $t^{-1}$ decay we perform an angular decomposition of the phase space into two regions.
One region contains the space resonant set $\mS$ but is away from the time resonant set $\mT$. 
In this region we can perform a normal form at the expense of introducing only a mild singularity when one of the Fourier variables vanishes.
For the quadratic boundary terms arising in the integration by parts in time, 
the decay is obtained in a straightforward manner, as it is done for the $g$ component. 
For the cubic terms the decay is obtained as a consequence of $L^2$-weighted estimates.
The complimentary region is away from $\mS \cap \mT^c$ and contains $\mR$. 
There we can combine the identity \eqref{xidxiphi} and the fact that $\phi$ can be divided by $\nabla_\eta \phi$, 
to conclude, roughly speaking, that $\nabla_\xi \phi \sim \nabla_\eta \phi$ in this region.
This implies a good control on weighted norms, and decay is obtained by interpolating in an appropriate fashion these norms.

The cubic terms $h_0(f,g)$ and $h_1$ are estimated in section \ref{weightedh_0fg},
using again the decomposition $f = g + h$.
Also for these terms the pointwise decay is a consequence of the $L^2$-weighted bounds. 
Some results from linear Harmonic analysis and the proof of Theorem \ref{theoop} are provided in the Appendix.

\subsection{Quantities controlled by the $X$ norm}
Here we give some estimates which will be useful in our proof. 
They follow from interpolating between the various components of the $X$-norm \eqref{Xnorm}.

\begin{lem}
For $0 < k < N-1$, and $2 \leq p \leq \infty$, one has
\begin{equation}
\label{dkuLp}
{\left\| \nabla^k u \right\|}_{L^p} \lesssim \frac{1}{ t^{\left(1-\frac{2}{p} \right)} }  t^{\alpha(k,N,p)} {\| u \|}_X
		\quad ,  \quad \mbox{where} \quad \alpha(k,N,p) := \frac{k \left( 1 - \frac{2}{p} + \ve \right) 
		}{N - \frac{3}{2} + \frac{3}{p}  } \, .
\end{equation}
In particular, for any $0 < k \leq 3$, and $\ve$ sufficiently small,
\begin{align}
\label{dNk}
\alpha(k,N,p) \leq 
		\frac{6}{N} := \d_{N} \, .
\end{align}
We also have
\begin{equation}
\label{uW14}
{\left\|  u \right\|}_{W^{1,p}} \lesssim \frac{1}{ t^{\left(1-\frac{2}{p} \right)} }  {\| u \|}_X 
\qquad \mbox{for $2 \leq p \leq 4$}  \, , 
\end{equation}
and
\begin{equation}
\label{u/gradLp}
{\left\|  \frac{u}{\grad} \right\|}_{L^p} \lesssim \frac{1}{ t^{\left(1-\frac{2}{p} \right)} } t^\g  {\| u \|}_X 
\qquad \mbox{for $4 \leq p < 6$} \, .
\end{equation}
\end{lem}

\proof
The proof of \eqref{dkuLp} follows from interpolation between the $H^N$ and the $L^\infty$ bounds
given by \eqref{Xnorm}. The proof is standard and can be found in \cite{GMS2}.
Inequality \eqref{uW14} follows from the dispersive estimate \eqref{dispp} for the linear wave propagator,
and from the definition of the $X$ norm \eqref{Xnorm}.
Similarly, it is easy to derive \eqref{u/gradLp} from \eqref{dispp} and the bounds on ${\| xf \|}_{H^1}$
provided by \eqref{Xnorm} $_\Box$ 

\section{Energy estimates}
\label{secenergy}
Energy estimates on $g$, $h_0$, and $h_1$ are straight forward.  To estimate $g$ we use 
Theorem \ref{theoop} and \eqref{fractional} to get
\begin{align*}
{\| g \|}_{H^N} & = {\left\| \grad^{-1} u \, u  \right\|}_{H^N}
\lesssim {\left\| \grad^{-1} e^{it\grad} f  \right\|}_{W^{1,9}}  {\| u \|}_{H^N}
+
{\left\| \grad^{-1} f \right\|}_{H^N}  {\| u \|}_{W^{1,6}}
\\
&
\lesssim {\left\| e^{it\grad} f \right\|}_{W^{1,\frac{9}{4}}}  {\| u \|}_{H^N}
+
\left( {\left\| f \right\|}_{H^{N-1}} +  {\left\| \grad^{-1} f \right\|}_{L^2} \right)  {\| u \|}_{W^{1,6}}
\\
&
\lesssim \frac{1}{t^\frac{1}{9}}  {\| u \|}_{X}  t^\ve {\| u \|}_{X} 
+ \left(  t^\ve {\| u \|}_{X} +  {\| xf \|}_{L^2} \right)  \frac{1}{t^\frac{2}{3}} t^{\d_{N}} {\| u \|}_{X}
\lesssim {\| u \|}_{X}^2 \, ,
\end{align*}
provided $\ve$ and $\g$ are small enough, and $N$ is large enough so that $\ve + \d_{N} < \frac{2}{3}$.
To estimate $h_0$ in $H^N$ we write
\begin{align*}
{\left\| h_0 \right\|}_{H^N} 
& = {\left\| \int_1^t  e^{-is\grad} T_{\nabla_\eta \phi(\xi,\eta)} \left( e^{-is\Lambda} f ,
				e^{is\Lambda} f \right) \, ds \right\|}_{H^N}
\\
& \lesssim  \int_1^t  {\left\|  e^{-is\Lambda} f \right\|}_{H^N}
			{\left\| e^{is\Lambda} R f \right\|}_{L^\infty}
 \lesssim {\| u \|}_{X}^2
	\int_1^t  s^\ve \frac{1}{s} \, ds
\lesssim
	t^\ve {\| u \|}_{X}^2 \, ,
\end{align*}
To estimate $h_0$ in $H^2$ we integrate by parts in $\eta$ and rewrite it by 
\begin{equation*}
\begin{split}
\what{h}_0 (t,\xi) &= \int_1^t \! \int \frac{1}{s}  e^{is\phi(\xi,\eta)} \nabla_\eta \what{f}(t,\eta) \what{f}(t,\xi-\eta)  \, d\eta ds
 \,\, + \,\, \mbox{symmetric term} \, ,
\end{split}
\end{equation*}
so to obtain
\begin{align*}
{\left\| h_0 \right\|}_{H^2} 
& = {\left\| \int_1^t \frac{1}{s} e^{-is\grad} \left( e^{-is\Lambda} x f \,
				e^{is\Lambda} f \right) \, ds \right\|}_{H^2}
\\
& \lesssim
	\int_1^t  \frac{1}{s} \left[ {\left\|  e^{-is\Lambda} x f \right\|}_{H^2}
			{\left\| e^{is\Lambda} f \right\|}_{W^{1,4}}  
	+
	{\left\|  e^{-is\Lambda} x f \right\|}_{L^6}
			{\left\| e^{is\Lambda} f \right\|}_{W^{2,3}}
	\right] \, ds
\\
& \lesssim
	{\| u \|}_{X}^2
	\int_1^t  \frac{1}{s} \left[\frac{s^\g}{\sqrt{s}}
	+
	\frac{1}{s^\frac{1}{3}} s^{\d_{N}} \right ]  \, ds
\lesssim {\| u \|}_{X}^2 \, ,
\end{align*}
provided $\d_{N} < \frac{1}{3}$.
Finally, to estimate $h_1$ note that
\begin{align*}
h_1 (t,x) 
& = \int_1^t e^{-is\grad} \left[  \left( \grad^{-1} e^{-is\grad} \partial_s \bar{f} \right) \, e^{is\grad} f
+
\left( \grad^{-1} e^{-is\grad} \bar{f} \right) \, e^{is\grad} \partial_s f  \right]
\, ds \, .
\end{align*}
Since $e^{is\grad} \partial_s f = T_q (f,f)$ for $q \in \B_0$,
we have that $e^{is\grad} \partial_s f \sim u^2$ as far as estimates are concerned.
One can then proceed in the same way as done for $g$ above to obtain
${\left\| h_0 \right\|}_{H^N} \lesssim {\| u \|}_X^3$.


\section{Weighted and $L^\infty$ estimates on $g$}\label{weightedg}

\subsection{Estimate of $\sup_t t^{-\g} {\| x g \|}_{L^2}$ and $\sup_t {\| x \grad g \|}_{H^1}$}
By Plancharel's Theorem estimating ${\| x g \| }_{H^2}$ is 
equivalent to estimate ${\| {\langle \xi \rangle}^2 \nabla_\xi \what{g} \| }_{L^2}$.
Applying $\nabla_\xi$ to $\what{g}$  we get:
\begin{align}
\nabla_\xi \what{g} (\xi) 
& = 
\label{dxig1}
\int t \, \nabla_\xi \phi(\xi,\eta) \, e^{it\phi} \frac{ \what{f}(\eta)}{|\eta|}  \what{f}(\xi-\eta) \, d\eta
+ \int e^{it\phi} \frac{\what{f}(\eta)}{|\eta|}  \nabla_\xi \what{f}(\xi-\eta) \, d\eta  = \hat I  + \widehat{II}\, .
\end{align}
%
%
%
Since $\nabla_\xi \phi \in \B_0$,
we can use Theorem \ref{theoop} and the dispersive estimate \eqref{dispp} to obtain
\begin{align*}
{\left\|I \right\|}_{L^2} 
	&  =  t {\left\| T_{\nabla_\xi \phi (\xi, \eta)} \left( e^{it\grad} \frac{f}{\grad} , u\right) \right\|}_{L^2}
 \lesssim  t {\left\| e^{it\Lambda} \frac{f}{\grad}  \right\|}_{L^4}  {\| u \|}_{L^4}
 \lesssim   \frac{t}{\sqrt{t}} {\| \langle x \rangle  f \|}_{L^2}  \frac{1}{\sqrt{t}} {\| u \|}_{X}
\lesssim  t^\g  {\| u \|}_{X}^2 \, .
\end{align*}
$II$  can be directly estimated in $H^2$ by using theorem \ref{theoop}, \eqref{fractional},
and Sobolev's embeddings:
\begin{align*}
{\left\|  II \right\|}_{H^2} 
	&	=  {\left\| e^{it\grad} \frac{f}{\grad}  \, e^{it\grad} x f \right\|}_{H^2}
 \lesssim  {\left\| e^{it\Lambda} \frac{f}{\grad} \right\|}_{W^{1,4}}  {\| e^{it\grad} x f \|}_{H^2}
		+ {\left\| e^{it\Lambda} \frac{f}{\grad}  \right\|}_{W^{2,4}}  {\| e^{it\grad} x f \|}_{L^3}		
\\
	& \lesssim \frac{1}{\sqrt{t}} {\| \langle x \rangle f \|}_{H^2}  {\| x f \|}_{H^2}
 \lesssim \frac{1}{\sqrt{t}} t^\g {\| u \|}_{X} t^\g {\| u \|}_{X} 
\lesssim {\| u \|}_{X}^2 \, .
\end{align*}

Since we have already estimated the $L^2$ norm of $\nabla_\xi \what{g} (\xi) $,
in order to estimate $\F^{-1}\nabla_\xi \what{g} (\xi) $ in $H^2$ we just need to bound $\F^{-1} |\xi| \nabla_\xi \what{g} (\xi) $ in $H^1$.
From \eqref{dxig1} we have
\begin{align*}
|\xi| \nabla_\xi \what{g} (\xi) 
& = 
\int e^{it\phi} t |\xi| \nabla_\xi \phi \frac{1}{|\eta|} \what{f}(\eta) \what{f}(\xi-\eta) \, d\eta
\, .
\end{align*}
%
%
%
Since $| \xi | \nabla_\xi \phi = \frac{\eta-\xi}{|\eta-\xi|} \phi + |\eta| \nabla_\eta \phi$,
integrate  by parts in $\eta$  to obtain
\begin{align}
|\xi| \nabla_\xi \what{g} (\xi) 
\label{xidxig1}
= &\int t e^{it\phi}  \frac{\phi}{|\eta|} \what{f}(\eta) \, \frac{\xi - \eta}{|\xi - \eta|} \what{f}(\xi-\eta) \, d\eta
 + \int e^{it\phi} \nabla_\eta \what{f}(\eta) \, \what{f}(\xi-\eta) \, d\eta  
\\&+  \mbox{``similar term''}  =  \hat I + \widehat{II} +  \mbox{``similar term''}\, ,\nn
\end{align}
where ``similar term'' denotes the term where $\nabla_\eta$ hits the other profile. 
Since $\frac{\phi}{|\eta|}$ is assumed to be in the class $\B_0$, we see that
$
I = t e^{-it\grad} T_{\mu_0 (\xi,\eta)} (u,u)
$,
with $\mu_0 \in \B_0$.
Then by Theorem \ref{theoop} and \eqref{uW14} we can estimate
\begin{align*}
{\left\| I\right\|}_{H^1} 
	&  =  t {\left\| T_{\mu_0(\xi,\eta)} (u, u) \right\|}_{H^1}
 \lesssim  t  {\| u \|}_{W^{1,4}}^2 
 \lesssim  {\| u \|}_{X}^2 \, .
\end{align*}
The term $II$ can be handle as follows:
\begin{align*}
{\left\| II\right\|}_{H^1} 
	&  =  {\left\| e^{it\grad} x f \, u \right\|}_{H^1}
 \lesssim  {\| e^{it\grad} x f \|}_{W^{1,6}}  {\| u \|}_{W^{1,3}}
 \lesssim  {\| \grad x f \|}_{H^1} {\| u \|}_{H^2}
\lesssim  {\| u \|}_{X}^2 \, .
\end{align*}

\subsection{Bounds on $ \sup_t t^{-1} {\| \grad x^2 g \|}_{H^1}$}


In the previous section we saw that $|\xi| \nabla_\xi \what{g}$ is made of the following
three types of contributions:
\[
\int t e^{it\phi} \mu_0(\xi,\eta) \what{f}(\eta) \what{f}(\xi-\eta) \, d\eta, \quad
 \int e^{it\phi} \nabla_\eta  \what{f}(\eta) \what{f}(\xi-\eta) \, d\eta, \quad
   \int |\xi|  e^{it\phi} \frac{1}{|\eta|}  \what{f}(\eta)  \nabla_\xi \what{f}(\xi-\eta) \, d\eta \, .
\]
In order to achieve the desired bound on $\grad x^2 g$, apply $\nabla_\xi$ to
the above terms and estimate the resulting expressions.
Using \eqref{xidxiphi} to deduce that $\frac{|\xi| \nabla_\xi \phi}{|\eta|} \in \B_0$,
we get contributions of the form
\begin{subequations}
\begin{align}
\label{xidxig11}
& \int t^2 e^{it\phi} \mu_0(\xi,\eta) \what{f}(\eta) \what{f}(\xi-\eta) \, d\eta
\\
\label{xidxig13}
& \int t e^{it\phi} \mu_0(\xi,\eta) \what{f}(\eta) \nabla_\xi  \what{f}(\xi-\eta) \, d\eta 
\\
\label{xidxig21}
&  \int e^{it\phi} \nabla_\eta  \what{f}(\eta) \nabla_\xi \what{f}(\xi-\eta) \, d\eta 
\\
\label{xidxig31}
&  \int |\xi| \, e^{it\phi} \frac{1}{|\eta|}  \what{f}(\eta)  \nabla_\xi^2 \what{f}(\xi-\eta) \, d\eta
\end{align}
\end{subequations}
plus similar or simper terms.  The first term is the one which gives the highest growth in $t$:
\begin{align*}
{\left\| \F^{-1} \eqref{xidxig11} \right\|}_{H^1} 
	&  =  t^2 {\left\| T_{\mu_0(\xi,\eta)} (u, u) \right\|}_{H^1}
 \lesssim  t^2  {\| u \|}_{W^{1,4}}^2
	\lesssim  t {\| u \|}_{X}^2 \, .
\end{align*}

We can take care of the contribution coming from  \eqref{xidxig13}  and  \eqref{xidxig21} with an $L^6 \times L^3$ estimate,
followed by Sobolev's embedding.  This will give
\[
{\left\| \F^{-1} \eqref{xidxig13} \right\|}_{H^1} +
{\left\| \F^{-1} \eqref{xidxig21} \right\|}_{H^1} \lesssim t{\| u \|}_{X}^2.
\]
Finally, always by means of theorem \ref{theoop}, Sobolev's embedding, and \eqref{fractional}, we can 
crudely estimate:
\begin{align*}
{\left\| \F^{-1} \eqref{xidxig31} \right\|}_{H^1} 
&  = {\left\| \grad T_{\mu_0(\xi,\eta)} \left( \grad^{-1} e^{it\grad} f,
				e^{it\grad} {|x|}^2 f \right) \right\|}_{H^1}
\\
& \lesssim   {\left\| \grad^{-1} e^{it\grad} f  \right\|}_{W^{2,6}} 
			{\left\| \grad {|x|}^2 f \right\|}_{H^1}
+ 
{\left\| e^{it\grad} f \right\|}_{W^{1,3}} 
			{\left\|{|x|}^2 f \right\|}_{W^{1,6}}
\\
& \lesssim {\left\| f \right\|}_{H^2}
		{\left\| \grad {|x|}^2 f \right\|}_{H^1}
\lesssim t {\| u \|}_{X}^2 \, .
\end{align*}

\subsection{$L^\infty$ estimate for $e^{i t \Lambda} g$\label{Linftyg}}
By  choosing  $N$ large enough 
\begin{align*}
{ \left\| e^{i t \Lambda} g  \right\| }_{\infty} 
& \lesssim {\left\| ( {\Lambda}^{-1} u) \, u  \right\| }_{W^{1,12}}
 \lesssim {\left\| {\Lambda}^{-1} u \right\|}_{W^{1,24}} {\| u \|}_{W^{1,24}}
\\
& \lesssim {\| u \|}_{L^\frac{8}{3} \cap L^{24}} {\| u \|}_{W^{1,24}}
 \lesssim  \frac{1}{ t^{\frac{1}{4}} } {\| u \|}_{X} \frac{1}{ t^\frac{11}{12} } t^{\d_{N}} {\| u \|}_{X}
	\lesssim \frac{1}{t} {\| u \|}^2_{X}.
\end{align*}

\section{Weighted estimates on  $h_0(h,h)$\label{secweightedh_0}}
Before estimating weighted norms of $h_0$ we need the following lemma:
\begin{lem}[Bounds on $\partial_t h$]\label{lem6}
Let $\d_{N}$ be the quantity defined in \eqref{dNk}. Then
\begin{align}
\label{bounddsf}
&{\left\|   \partial_t h  \right\|}_{H^2} \lesssim \frac{1}{t^{2-}} t^{\g + \d_{N}} {\| u \|}_{X}^2,
\\
\label{boundxdsf}
&{\left\|  x \partial_t h  \right\|}_{H^2} \lesssim \frac{1}{t^{1-}} t^{\g + \d_{N}} {\| u \|}_{X}^2 \, ,
\end{align}
where $c -$ denotes a number smaller but arbitrarily close to $c$.
\end{lem}

\proof
Recall that $h = h_0 + h_1$ with $h_0$ and $h_1$ defined by \eqref{h0hat} and \eqref{h1hat} respectively.
To estimate $\partial_th_0$  let $\d$ be arbitrarily small and estimate
\begin{align*}
{\left\| \partial_t h_0 \right\|}_{H^2} 
& = {\frac{1}{it}  \left\| e^{-it\grad} \left( e^{-it\Lambda} x f \,
				e^{it\Lambda} f \right) \right\|}_{H^2}
 \lesssim
	\frac{1}{t} {\left\|  e^{-it\Lambda} x f \right\|}_{H^2}
			{\left\| e^{it\Lambda} f \right\|}_{W^{3,\frac{1}{\d}}}
\\
& \lesssim
	\frac{1}{t} t^\g {\| u \|}_{X} \frac{1}{t^{1-2\d}}  t^{\d_{N}} {\| u \|}_{X}
 =
	\frac{1}{t^{2-}} t^{\g + \d_{N}} {\| u \|}_{X}^2 \, .
\end{align*}

To prove \eqref{boundxdsf} for $h_0$, we apply $\nabla_\xi$ to $\partial_t h_0$,
obtaining the following terms:
\begin{subequations}
\begin{align}
&
\label{dsh_03}
\int e^{it\phi(\xi,\eta)}  \nabla_\eta \phi(\xi,\eta) \what{f}(t,\eta) \nabla_\xi  \what{f}(t,\xi-\eta)  \, d\eta
\\
&
\label{dsh_02}
\int e^{it\phi(\xi,\eta)}  \nabla_\xi \nabla_\eta \phi(\xi,\eta)  \what{f}(t,\eta) \what{f}(t,\xi-\eta)  \, d\eta
\\
\label{dsh_01}
& \int it \nabla_\xi \phi e^{it\phi(\xi,\eta)} \nabla_\eta \phi(\xi,\eta) \what{f}(t,\eta) \what{f}(t,\xi-\eta)  \, d\eta .
\end{align}
\end{subequations}
To bound  \eqref{dsh_03}, we note that it is of the type
\begin{equation*}
\int e^{it\phi(\xi,\eta)}  \mu_0 (\xi,\eta) \what{f}(t,\eta) \nabla_\xi  \what{f}(t,\xi-\eta)  \, d\eta.
\end{equation*}
So let $\d$ be an arbitrarily small number and use Theorem \eqref{theoop} to get
\begin{align*}
{\left\| \eqref{dsh_03} \right\|}_{H^2} 
& = {\left\| e^{-it\grad} T_{\mu_0 (\xi,\eta)} \left( e^{-it\Lambda} f ,
				e^{it\Lambda} x f \right) \right\|}_{H^2}
 \lesssim
	{\left\|  e^{-it\Lambda} f \right\|}_{ W^{3,\frac{1}{\d}} }
			{\left\|  x f \right\|}_{H^2}
\\
& \lesssim
	\frac{1}{t^{1-2\d}}  t^{\d_{N}} {\| u \|}_{X}   t^{\g} {\| u \|}_{X}
=
	\frac{1}{t^{1-}} t^{\g + \d_{N}} {\| u \|}_{X}^2 \, .
\end{align*}
The term \eqref{dsh_02} is of the form
\begin{equation*}
\int e^{it\phi(\xi,\eta)}  \what{f}(t,\eta) \frac{\widehat{R^2 f}(t,\xi-\eta)}{|\xi-\eta|}  \, d\eta,
\end{equation*}
so by Hardy's inequality it can be estimated in the same fashion as \eqref{dsh_03}.
To estimate  \eqref{dsh_01}, we integrating by parts to get  terms like \eqref{dsh_02} and \eqref{dsh_03}.

Since $e^{it\grad} \partial_t \what{f}$ is a quadratic expression in $\what{f}$,
$h_1$ is essentially a cubic term (with a singularity of the type $\grad^{-1}$ on one of the three profiles).
The presence of an extra term is basically equivalent to the gain of $t^{-1}$ present in $\partial_t h_0$.
Therefore, the bounds \eqref{bounddsf} and \eqref{boundxdsf} for $h_1$ can be proven similarly 
to what we did above for $h_0$. We skip the details $_\Box$

\subsection{Estimate of ${ \sup_t {\| x h_0 \| }_{H^2}}$\label{secxh_0}}
For simplicity of notations, in the remainder of this section we will denote $h_0 (h,h)$ simply by $h_0$.   Recall that
\begin{equation*}
\begin{split}
\what{h}_0 (t,\xi) & \stackrel{def}{=} \int_1^t \! \int e^{is\phi(\xi,\eta)} 
	\nabla_\eta \phi(\xi,\eta)  \what{h}(s,\xi-\eta) \what{h}(s,\eta) \, d\eta ds,
\\
& = 
\int_1^t \! \int \frac{1}{s}  e^{is\phi(\xi,\eta)} \nabla_\eta \what{h}(t,\eta) \what{h}(t,\xi-\eta)  \, d\eta ds
 \quad + \quad \mbox{symmetric term} \, .
\end{split}
\end{equation*}By Plancharel's Theorem, estimating ${\| x h_0 \| }_{H^2}$ is 
equivalent to estimate ${\| {\langle \xi \rangle}^2 \nabla_\xi \what{h}_0 \| }_{L^2}$.
To this end, we first compute, using \eqref{xidxiphi}, an expression for $|\xi| \nabla_\xi \what{h}_0$:
\begin{align}
\label{xidxids}
|\xi| \nabla_\xi \what{h}_0 (\xi) 
& = \int_1^t \! \int s \, \partial_s e^{is \phi}  \nabla_\eta \phi(\xi,\eta) 
					\what{h} (\eta) \widehat {R h} (\xi-\eta)  \, d\eta \, ds
\\
\label{xidxideta}
& + \int_1^t \! \int s \, |\eta| \nabla_\eta e^{is \phi}  \nabla_\eta \phi(\xi,\eta) 
					\what{h} (\eta) \what{h} (\xi-\eta)  \, d\eta \, ds
\\
\nn
& + \int_1^t \! \int e^{is \phi}  |\xi| \nabla_\xi  \nabla_\eta \phi(\xi,\eta)
					\what{h} (\eta) \what{h} (\xi-\eta)  \, d\eta \, ds
\\
\nn
& + \int_1^t \! \int  e^{is \phi}  \nabla_\eta \phi(\xi,\eta)
			\what{h} (\eta) |\xi| \nabla_\xi \what{h} (\xi-\eta)  \, d\eta \, ds \, .
\end{align}
Integrating by parts in time in \eqref{xidxids},
and in frequency in \eqref{xidxideta}, we obtain (after collecting terms appropriately):
\begin{subequations}\label{g6}
\begin{align}
\label{xidxi1}
|\xi| \nabla_\xi \what{h}_0 (\xi) 
& \sim \int \! t  \, e^{it \phi}  \nabla_\eta \phi(\xi,\eta) 
						\what{h} (\eta) \widehat{R h} (\xi-\eta)  \, d\eta
\\
\label{xidxi2}
& -  \int_1^t \! \int s \, e^{is \phi} \nabla_\eta \phi(\xi,\eta)
				\partial_s \what{h} (\eta) \widehat{R h} (\xi-\eta)  \, d\eta \, ds
\\
\label{xidxi3}
& + \int_1^t \! \int e^{is \phi}  \left[ |\xi| \nabla_\xi - |\eta| \nabla_\eta \right] \nabla_\eta \phi(\xi,\eta)
					\what{h} (\eta) \what{h} (\xi-\eta)  \, d\eta \, ds
\\
\label{xidxi4}
& - \int_1^t \! \int e^{is \phi}  \nabla_\eta \phi(\xi,\eta)
			|\eta| \nabla_\eta \what{h} (\eta) \what{h} (\xi-\eta)  \, d\eta \, ds \, .
\\
\label{xidxi5}
& + \int_1^t \! \int e^{is \phi}   \nabla_\eta \phi(\xi,\eta) \what{h} (\eta) 
			\left[ |\xi| \nabla_\xi - |\eta| \nabla_\eta \right]  \what{h} (\xi-\eta)  \, d\eta \, ds \, 
\end{align}
\end{subequations}
plus   a boundary integral( \eqref{xidxi1} at $t=1$),   
a term symmetric to \eqref{xidxi2} corresponding to the case where $\partial_s$ hits the other profile,
and a term obtained when $\nabla_\eta$ hits $|\eta|$. The bound on these  terms are either similar or  easier than \eqref{g6}  and thus will be ignored.

\nl
We now make some observations that will further simplify our calculations:

\vskip 5pt

\nl (1)  In \eqref{xidxi5} we can write
\begin{equation*}
\left[ |\xi| \nabla_\xi - |\eta| \nabla_\eta \right]  \what{f} (\xi-\eta)
= |\xi - \eta| \nabla_\xi \what{f} (\xi-\eta)
	- \phi(\xi,\eta) \nabla_\xi \what{f} (\xi-\eta) \, .
\end{equation*}
The integral corresponding to the first summand above 
will give a term analogous to \eqref{xidxi4},
and can therefore be considered among the ``similar terms''.
The second summand can instead be treated by integration by parts in time,
yielding ``easier terms''.
\vskip 5pt
\nl (2)   Using \eqref{xidxiphi}, one can see that
\begin{equation*}
\left[ | \xi | \nabla_\xi - |\eta| \nabla_\eta  \right]  \nabla_\eta \phi (\xi, \eta) 
	=  \mu_0 (\xi,\eta) \nabla_{\eta} \phi + \mu_0^\p (\xi,\eta) \frac{\phi (\xi,\eta) }{|\xi-\eta|} \, .
\end{equation*}
Using the above identity we see that
\begin{align}
\label{xidxi3.1}
\eqref{xidxi3}	
& = \int_1^t \! \int e^{is \phi}  \mu_0(\xi,\eta) \nabla_{\eta} \phi (\xi, \eta) 
				\what{h} (\eta) \what{h} (\xi-\eta)  \, d\eta \, ds
\\
\label{xidxi3.2}
& + \left. \int e^{is \phi}  \mu_0^\p (\xi,\eta) \frac{1}{|\xi-\eta|}
			\what{h} (\eta) \what{h} (\xi-\eta)  \, d\eta \, ds \right]_1^t
\\
\label{xidxi3.3}
& - \int_1^t \! \int e^{is \phi}  \mu_0^\p (\xi,\eta) \frac{1}{|\xi-\eta|}
			\partial_s \left( \what{h} (\eta) \what{h} (\xi-\eta)  \right) \, d\eta \, ds \, .
\end{align}
All of these terms can be considered among the ``similar and easier terms'':
\eqref{xidxi3.1} is essentially like $\what{h}_0 (h,h)$,
while \eqref{xidxi3.2}, respectively \eqref{xidxi3.3},
is similar to one of the contributions that will appear when integrating by parts in frequency 
in \eqref{xidxi1}, respectively in \eqref{xidxi2}.
In view of the above observations, we write
\begin{equation}
\label{xidxiconv}
|\xi| \nabla_\xi \what{h}_0 (t,\xi) \sim \eqref{xidxi1} + \eqref{xidxi2} + \eqref{xidxi4}  =: \hat A(t,\xi) \, .
\end{equation}

\paragraph{\it Estimate of ${\| x h_0 \|}_{L^2}$}  From \eqref{xidxiconv} we have 
\begin{align}
\nonumber 
{\| x h_0 \|}_{L^2}& = {\left\|  \F^{-1} \nabla_\xi \what{h}_0 \right\|}_{L^2} 
	= {\left\| {\Lambda} ^{-1}\F^{-1} |\xi| \nabla_\xi \what{h}_0 \right\|}_{L^2}
	\sim {\left\| {\Lambda}^{-1}A\right\|}_{L^2} \, ,
\end{align}
and since all the bilinear terms present have  $\nabla_\eta \phi$ in their symbol, we can integrate by parts in frequency before we estimate. Thus 
\begin{align*}
e^{it\Lambda} \F^{-1} \eqref{xidxi1} 
			=  T_{\mu_0 (\xi, \eta)} \left(e^{-it\Lambda} x h, R e^{it\Lambda} h \right) 
			+  T_{\mu_0 (\xi, \eta)} \left(e^{-it\Lambda} h, e^{it\Lambda}  x R h \right) \, ,
\end{align*}
and using \eqref{fractional2} we can estimate
\begin{align*}
{\left\| \frac{1}{\Lambda}\F^{-1} \eqref{xidxi1}  \right\|}_{L^2} 
	& \lesssim   {\left\| T_{\mu_0 (\xi, \eta)} (e^{-it\Lambda} x h, R e^{it\Lambda} h) \right\|}_{L^\frac{6}{5}}
		+ {\left\| T_{\mu_0 (\xi, \eta)} ( e^{-it\Lambda} h, e^{it\Lambda}  x R h) \right\|}_{L^\frac{6}{5}}
\\
	& \lesssim   {\left\| x h \right\|}_{L^2}  {\| e^{it\Lambda} h \|}_{L^3}
		+ {\| e^{it\Lambda} h \|}_{L^3}  {\left\| e^{it\Lambda} x R h \right\|}_{L^2} 
 \lesssim  t^\g {\| u \|}_{X} \frac{1}{t^{1/3}}  {\| u \|}_{X} \lesssim {\| u \|}_{X}^2 \, .
\end{align*}

Similarly, by using lemma \ref{lem6} 
we get
\begin{align*}
{\left\| \frac{1}{\Lambda}\F^{-1} \eqref{xidxi2} \right\|}_{L^2} 
	& \lesssim  \int_1^t  {\left\| T_{\mu_0 (\xi, \eta)} \left(e^{-is\Lambda} x \partial_s h, 
				R e^{is\Lambda} h \right) \right\|}_{L^\frac{6}{5}} \, ds
\\
	& + \int_1^t  {\left\| T_{\mu_0 (\xi, \eta)} \left(e^{-is\Lambda} \partial_s h,
				e^{is\Lambda} x R h \right) \right\|}_{L^\frac{6}{5}} \, ds
\\
	& \lesssim  \int_1^t  {\left\| x \partial_s h \right\|}_{L^2}  {\| e^{is\Lambda} h \|}_{L^3} \, ds
		+  \int_1^t  {\left\| e^{-is\Lambda} \partial_s h \right\|}_{L^3}
			{\left\| x R h \right\|}_{L^2} \, ds
\\
	& \lesssim 
{\| u \|}_{X}^3 \int_1^t \frac{1}{s^{1-}}  s^{\d_{N} + \g}  \frac{1}{s^\frac{1}{3}} \, ds  
	\lesssim {\| u \|}_{X}^3  \, ,
	\\
%
{\left\| \frac{1}{\Lambda}\F^{-1} \eqref{xidxi4} \right\|}_{L^2} 
	& \lesssim  \int_1^t  \frac{1}{s} {\left\| T_{\mu_0 (\xi, \eta)} \left(e^{-is\Lambda} \Lambda x h,
			e^{is\Lambda} x h \right) \right\|}_{L^\frac{6}{5}} \, ds
\\
	& + \int_1^t  \frac{1}{s} {\left\| T_{\mu_0 (\xi, \eta)} \left(e^{-is\Lambda} \Lambda {|x|}^2 h,  
			e^{is\Lambda} h \right) \right\|}_{L^\frac{6}{5}} \, ds
\\
	& \lesssim    \int_1^t  \frac{1}{s}  {\left\| e^{-is\Lambda} \Lambda x h \right\|}_{L^3}  {\| x h \|}_{L^2} \, ds
		+
	\int_1^t  \frac{1}{s} {\left\| {\Lambda}{|x|}^2 h \right\|}_{L^2}  {\left\| e^{is\Lambda} h \right\|}_{L^3} \, ds
\\
	& \lesssim 
		{\| u \|}_{X}^2 \int_1^t \frac{1}{s} \frac{1}{s^{1/3}} s^a  s^\g  \, ds
		+
		{\| u \|}_{X}^2 \int_1^t \frac{1}{s} \frac{1}{s^{1/3}} s^a  \, ds
	\lesssim {\| u \|}_{X}^2 \, .
\end{align*}

\paragraph{ \it Estimate of ${\| \Lambda x h_0 \|}_{H^1}$} This  is
similar  to estimating 
$\F^{-1} $ of \eqref{xidxi1}, \eqref{xidxi2}, and \eqref{xidxi4}   in $H^1$.
This can be done in a similar fashion to the previous paragraph.  In fact these
estimates  are even easier,  since the ${\Lambda}^{-1}$ singularity is not present any more.
We only show how to bound ${\| \F^{-1} \eqref{xidxi4} \|}_{H^1}$, the other estimates being very similar.  By lemma \ref{lem6}, we have
\begin{align*}
& {\left\|  \F^{-1} \eqref{xidxi4}\right\|}_{H^1}
\\
	& \lesssim  
	\int_1^t  {\left\| T_{\mu_0 (\xi, \eta)} \left(e^{-is\Lambda} x \partial_s h, 
						R e^{is\Lambda} h \right) \right\|}_{H^1} ds
+ \int_1^t  {\left\| T_{\mu_0 (\xi, \eta)} \left(e^{-is\Lambda} \partial_s h, 
		e^{is\Lambda} x R h \right) \right\|}_{H^1} ds
\\
	& \lesssim  \int_1^t  {\left\| x \partial_s h \right\|}_{W^{1,6}}  
		{\| e^{is\Lambda} h \|}_{W^{1,3}} ds
	+
		\int_1^t  {\left\| e^{-is\Lambda} \partial_s h \right\|}_{W^{1,6}} 
		{\left\| x R h \right\|}_{H^1} ds
\\
	& \lesssim  {\| u \|}_{X}^3 \int_1^s \frac{1}{s^{1-}} s^{\g+\delta_{N}} 
		\frac{1}{s^\frac{1}{3}}  \, ds
		\lesssim {\| u \|}_{X}^3  \, .
\end{align*}

\subsection{Estimate of $\sup_t t^{-a} {\left\| {|x|}^2 \grad h_0 \right\| }_{L^2}$}  
By reducing the bilinear interacting to the prototypes given in   \eqref{xidxi1}, \eqref{xidxi2}, and \eqref{xidxi4}.  we  proceed as follows
Applying $\nabla_\xi$ to \eqref{xidxi1} produces terms of the form:
\begin{subequations}\label{dxi615}
\begin{align}
\label{dxixidxi11}
& \int \! t^2  e^{it \phi (\xi,\eta)} \, \mu_0 (\xi,\eta) \, \nabla_\eta \phi(\xi,\eta)
				\what{h} (\eta) \what{h} (\xi-\eta)  \, d\eta
\\
\label{dxixidxi12}
&  \int \! t \, e^{it \phi (\xi,\eta)} \, \mu_0 (\xi,\eta) \, \nabla_\xi \nabla_\eta \phi(\xi,\eta)
		\what{h} (\eta) \what{h} (\xi-\eta)  \, d\eta
\\
\label{dxixidxi13}
& \int \! t \, e^{it \phi (\xi,\eta)} \, \mu_0 (\xi,\eta) \, \nabla_\eta \phi(\xi,\eta)
		\what{h} (\eta) \nabla_\xi  \what{h} (\xi-\eta)  \, d\eta  \, .
\end{align}
\end{subequations}
Integrating by parts in $\eta$ whenever there is  $\nabla_\eta \phi$ reduces \eqref{dxi615} to terms of the following type
\begin{subequations}\label{dxixidxi1110}

\begin{align}
\label{dxixidxi111}
& \int \! t   \, e^{it \phi (\xi,\eta)}  \nabla_\eta \mu_0(\xi, \eta) \what{h} (\eta) \what{h} (\xi-\eta)   \, d\eta
\\
\label{dxixidxi112}
& \int \! t  \, e^{it \phi(\xi,\eta)}  \mu_0(\xi, \eta) \nabla_\eta \what{h} (\eta) \what{h} (\xi-\eta)  \, d\eta \, .
\end{align}
\end{subequations}
To estimate \eqref{dxixidxi111} we note that 
 for any symbol $\mu_0 \in \B_0$,  $ \nabla_\eta \mu_0 (\xi,\eta) = \frac{\mu_0^{\p}(\xi,\eta)}{|\xi-\eta|} + \frac{\mu_0^{\p\p}(\xi,\eta)}{|\eta|} $
for some symbols $\mu_0^\p , \mu_0^{\p\p} \in \B_0$, consequently  
\begin{align*}
{\| \eqref{dxixidxi111} \|}_{L^2}
& \lesssim
	t {\left\| e^{-it\Lambda} h \right\|}_{L^4}
	{\left\|e^{it\Lambda} \grad^{-1} h \right\|}_{L^4}
 \lesssim
	t \frac{1}{\sqrt{t}} {\left\| u \right\|}_X  
		\frac{1}{\sqrt{t}} {\left\| \langle x \rangle h \right\|}_{ L^2 }
 \lesssim
		t^\g {\left\| u \right\|}^2_X
	\lesssim  t^{a} {\left\| u \right\|}^2_X 
\end{align*}
since  $ \g \leq a$.
We can deal similarly with \eqref{dxixidxi112}:
\begin{align*}
{\| \eqref{dxixidxi112} \|}_{L^2} 
& \lesssim
		t {\left\| e^{-it\Lambda} h \right\|}_{L^4}
 \lesssim  
		t \frac{1}{\sqrt{t}} {\left\| u \right\|}_X \frac{1}{\sqrt{t}} {\left\| \grad {|x|}^2 h \right\|}_{L^2}
		\lesssim  t^{a} {\left\| u \right\|}^2_X   \, .
\end{align*}

Next  apply $\nabla_\xi$ to \eqref{xidxi2}  and integrate by parts in $\eta$ whenever   $\nabla_\eta \phi$ is present in the symbol 
to get  terms of the  type:
\begin{subequations}\label{dxixidxi240}
\begin{align}
\label{dxixidxi24}
& \int_1^t \! \int s \, e^{is \phi(\xi,\eta)}  \mu_0 (\xi, \eta) 
		\partial_s  \nabla_\eta \what{h} (\eta) \what{h} (\xi-\eta)  \, d\eta \, ds
\\
\label{dxixidxi22}
&  \int_1^t \! \int s \, e^{is \phi(\xi,\eta)}  \mu_0 (\xi, \eta) 
		\partial_s \what{h} (\eta) \frac{\what{h} (\xi-\eta)}{|\xi-\eta|}  \, d\eta \, ds
\\
\label{dxixidxi23}
&  \int_1^t \! \int s \, e^{is \phi(\xi,\eta)}  \mu_0 (\xi, \eta) 
		\partial_s \what{h} (\eta) \nabla_\xi \what{h} (\xi-\eta)  \, d\eta \, ds
\end{align}
\end{subequations}
From lemma \ref{lem6}, and choosing a small enough $\d$, we have
\begin{align*}
{\left\| \eqref{dxixidxi24} \right\|}_{L^2} 
& \lesssim   \int_1^t s \,  {\left\|
	T_{ \mu_0 (\xi, \eta)} \left(e^{-is\Lambda} x \partial_s h, e^{is\Lambda} h  \right) \right\|}_{L^2}ds
 \lesssim  
	\int_1^t  s \, {\left\| e^{-is\Lambda} x \partial_s h \right\|}_{ L^{\frac{2}{1-2\d }}}
		{\| e^{is\Lambda} h \|}_{L^\frac{1}{\d}} ds
\\
& \lesssim 
	\int_1^t  s \, {\left\| x \partial_s h \right\|}_{H^1}
		\frac{1}{s^{1-2\delta}} {\| u \|}_X \, ds
 \lesssim			
	{\| u \|}_{X}^3 \int_1^t s \frac{1}{s^{1-}} s^{\g + \d_{N}} \frac{1}{s^{1-2\delta}} \, ds
\lesssim
	{\| u \|}_{X}^3 t^a \, ,\\
	{\left\| \eqref{dxixidxi22} \right\|}_{L^2} 
& \lesssim   \int_1^t s \, {\left\|
	T_{ \mu_0 (\xi, \eta)} \left(e^{-is\Lambda} \partial_s h,
		e^{is\Lambda} \frac{h}{\Lambda} \right) \right\|}_{L^2}   \, ds
 \lesssim  
	\int_1^t  s \, {\left\| e^{-is\Lambda} \partial_s h \right\|}_{L^\frac{1}{\d}} 
		{\left\| e^{is\Lambda} \frac{h}{\Lambda} \right\|}_{{ L^{\frac{2}{1-2\d }}}} \, ds
\\
& \lesssim
	\int_1^t  s  \,  \frac{1}{s^{2 (1 - \d)}} {\| u \|}_{X}^2
		{\| \langle x \rangle h\|}_2 \, ds
\lesssim {\| u \|}_{X}^3 t^a
\end{align*}
provided we choose $N$ large enough so that $\d_{N} + \g < a$.
Observe that \eqref{dxixidxi22} and \eqref{dxixidxi23} are similar to each other
through Hardy's inequality.

Finally to estimate  $\nabla_\xi \eqref{xidxi4}$ we have terms of the type
\begin{subequations}
\begin{align}
& \int_1^t \! \int s \nabla_\xi \phi \, e^{is \phi}  \nabla_\eta \phi \,
		|\eta| \nabla_\eta \what{h} (\eta) \what{h} (\xi-\eta)  \, d\eta \, ds
\label{dxixidxi41}
\\
\label{dxixidxi42}
& \int_1^t \! \int e^{is \phi}  \nabla_\xi \nabla_\eta \phi \,
		|\eta| \nabla_\eta \what{h} (\eta) \what{h} (\xi-\eta)  \, d\eta \, ds
\\
\label{dxixidxi43}
& \int_1^t \! \int e^{is \phi}   \nabla_\eta \phi \,
		|\eta| \nabla_\eta \what{h} (\eta) 
		\nabla_\xi \what{h} (\xi-\eta)  \, d\eta \, ds
\end{align}
\end{subequations}
Integrating by parts  in $\eta$ leads to the following  types of bilinear forms
\begin{subequations}
\begin{align}
& \int_1^t \! \int \nabla_\xi \phi \, e^{is \phi}  \,
		|\eta| \nabla_\eta^2 \what{h} (\eta) \what{h} (\xi-\eta)  \, d\eta \, ds
\label{dxixidxi411}
\\
\label{dxixidxi412}
& \int_1^t \! \int \nabla_\xi \phi \, e^{is \phi}  \,
		|\eta| \nabla_\eta \what{h} (\eta)  \nabla_\eta \what{h} (\xi-\eta)  \, d\eta \, ds
\end{align}
\end{subequations}
We can  bound \eqref{dxixidxi411} by an $L^2 \times L^\infty$ estimate,
since we have control over $R h$ in $L^\infty$ with a sharp bound of $t^{-1}$, see \eqref{RkuLinfty}. 
To estimate \eqref{dxixidxi412}
we interpolate between the available weighted $L^2$ bounds in the $X$ norm \eqref{boundsh}
to obtain
\begin{align*}
{\left\| \eqref{dxixidxi412} \right\|}_{L^2} 
& \lesssim \int_1^t {\left\| T_{ \mu_0 (\xi, \eta)} \left(e^{-is\Lambda} \Lambda {|x|} h,
				e^{is\Lambda} x h \right) \right\|}_{L^2} \, ds
 \lesssim
	\int_1^t   {\left\|  e^{-is\Lambda} \Lambda |x| h \right\|}_{L^4}
			{\left\| e^{is\Lambda} x h \right\|}_{L^4}  \, ds
\\
& \lesssim
 		\int_1^t  \frac{1}{\sqrt{s}}	{\left\| |x| h \right\|}_{\dot{W}^{2,\frac{4}{3}} }
 			\frac{1}{\sqrt{s}}  {\left\| x h \right\|}_{\dot{W}^{1,\frac{4}{3}} }  \, ds
 \lesssim
 		\int_1^t  \frac{1}{s}	{\left\| {|x|}^{\frac{7}{4} +} h \right\|}_{\dot{H}^2}
 			{\left\|  {\langle x \rangle}^{\frac{7}{4} +} h \right\|}_{\dot{H}^1}  \, ds
\\
& \lesssim
			{\| u \|}_{X}^2 \int_1^t \frac{1}{s} s^{\frac{3}{4} b +} 
			s^{\frac{3}{4} a +}  \, ds
\lesssim 
			{\| u \|}_{X}^2 t^a \, .
\end{align*}

\subsection{Estimate of $\sup_t t^{-b} {\left\| {|x|}^2 {\Lambda}^2 h_0 \right\| }_{L^2}$}
In order to bound ${|x|}^2 {\Lambda}^2 h_0$ in $L^2$ we need to estimate the $L^2$ norm of
$\xi\nabla_\xi$ of\ \eqref{xidxi1},  \ \eqref{xidxi2} and \eqref{xidxi4}.

To estimate  ${ \left\| \xi\nabla_\xi \eqref{xidxi1} \right\| }_{L^2}$ and
 ${ \left\| \xi\nabla_\xi \eqref{xidxi2} \right\| }_{L^2}$  we 
recall that  all the terms  reduce to  \eqref{dxixidxi1110},  and  \eqref{dxixidxi240}.  For these terms 
it is possible to proceed in the exact same way as we did before in the previous paragraph,
since the presence of an extra derivative does not cause any harm.
We just show how to deal with \eqref{dxixidxi111}:
\begin{align*}
{\left\|\xi \eqref{dxixidxi111} \right\|}_{L^2}
& \lesssim
	t {\left\| e^{-it\Lambda} h \right\|}_{ W^{1,\frac{1}{\d}} }
	{\left\|e^{it\Lambda} \grad^{-1} h \right\|}_{ W^{1, \frac{2}{1-2\d} } }
\\
& \lesssim
	t {\left\| u \right\|}_X \frac{1}{t^{1-2\delta}} t^{\d_{N}} 
		\left[ {\left\| h \right\|}_{ L^{ \frac{2}{1-2\d} } }
		+ {\left\| h \right\|}_{ L^{ \frac{5}{5-6\d} } } \right]
\\
	& \lesssim
		t^{2\delta + \d_{N}} {\left\| u \right\|}_X
		{\left\| \langle x \rangle h \right\|}_{H^1}
 \lesssim
		t^{2\delta + \d_{N}} {\left\| u \right\|}_X    t^\g {\left\| u \right\|}_X
\lesssim  t^{a} {\left\| u \right\|}^2_X 
\end{align*}
since we can choose $N$ large enough so that $\g + \d_{N} \le a$.

To estimate t \eqref{xidxi4}   apply $|\xi| \nabla_\xi$ to the equation, using  \eqref{xidxiphi}, to get 
\begin{subequations}
\begin{align}
\label{dxixidxi50}
&  \int  t \, e^{i t \phi} \,  \nabla_\eta \phi \,
		|\eta|  \nabla_\eta \what{h} (\eta) \what{R h} (\xi-\eta)  \, d\eta \, ds
\\
\label{dxixidxi51}
& \int_1^t \! \int  e^{is \phi}  \nabla_\eta \phi \,
		{|\eta|}^2  \nabla^2_\eta \what{h} (\eta) \what{h} (\xi-\eta)  \, d\eta \, ds
\\
\label{dxixidxi52}
& \int_1^t \! \int e^{is \phi} \nabla_\eta \phi \,
		|\eta| \nabla_\eta \what{h} (\eta) 
		\, |\xi-\eta| \nabla_\eta \what{h} (\xi-\eta)  \, d\eta \, ds
\\
\label{dxixidxi53}
& \int_1^t \! \int s \, e^{is \phi}  \nabla_\eta \phi \,
		|\eta| \nabla_\eta \partial_s \what{h} (\eta) 
		\what{R h} (\xi-\eta)  \, d\eta \, ds
\\
\nn
& + \quad \mbox{``similar and easier terms''} \, .
\end{align}
\end{subequations}
The term \eqref{dxixidxi50}, respectively \eqref{dxixidxi51},
can be bounded by $O(t^b)$,
using respectively an $L^4 \times L^4$ and an  $L^2 \times L^\infty$ estimate
(recall that we have control of $R h$ in $L^\infty$).

To estimate \eqref{dxixidxi52}, we use the presence of $\nabla_\eta \phi$ to integrate by parts.
Up to ``similar or easier terms'', this gives
\begin{equation}
\label{dxixidxi521}
\int_1^t \! \int \frac{1}{s} e^{is \phi}
		|\eta| \nabla_\eta^2 \what{h} (\eta) \, |\xi-\eta| \nabla_\eta \what{h} (\xi-\eta)  \, d\eta \, ds \, .
\end{equation}
Using Sobolev's embeddings, we can estimate
\begin{align*}
{\left\| \eqref{dxixidxi521} \right\|}_{L^2} 
& \lesssim 
		\int_1^t  \frac{1}{s}  {\left\|  e^{-is\Lambda} \Lambda {|x|}^2 f \right\|}_{L^6}
		{\left\|  e^{is\Lambda} \Lambda x  f \right\|}_{L^3}  \, ds
 \lesssim 
		\int_1^t  \frac{1}{s}  {\left\|  {\Lambda}^2 {|x|}^2 f \right\|}_{L^2}
		{\left\|  \Lambda  x  f \right\|}_{H^1}  \, ds
\\
& \lesssim	{\left\|  u \right\|}_X^2	 
		\int_1^t  \frac{1}{s} s^b \, ds \lesssim  {\left\|  u \right\|}_X^2	t^b \, . 
\end{align*}

Finally, observe that \eqref{dxixidxi53} is similar to \eqref{dxixidxi24}
and can be estimated in an analogous fashion. 
This concludes the proof of weighted estimates on $h_0$.



\section{$L^\infty$ estimate for $e^{i t \grad} h_0(h,h)$\label{Linftyh}}

Although the term $h_0$ is tight and satisfies stronger weighted estimates than $g$, these bounds are not good enough  
 to imply that $e^{i t \Lambda} h_0(h,h)$ decays like $t^{-1}$ in $L^\infty$.
We will achieve such a  bound  by dividing the interactions according to the resonance analysis carried out in section  
\ref{secresonances}.

\subsection{Angular partition of the phase space}
Let us introduce  the cutoff function 
\begin{equation}
\label{angcutoff}
\chi (\xi, \eta) := \wt{\chi} \left(  \frac{\xi}{|\xi|} \cdot \frac{\eta - \xi}{|\eta - \xi|}  \right)
\end{equation}
with $\wt{\chi}$ smooth non--decreasing function such that
\begin{equation}
\wt{\chi} (x) 
= \left\{
\begin{array}{ll}
1 & \mbox{for} \quad x \geq \frac{1}{4}
\\
0 & \mbox{for} \quad x \leq - \frac{1}{4} \, .
\end{array}
\right.
\end{equation}
Observe  that on the support of $\chi$  
we have $\frac{\xi}{|\xi|} \cdot \frac{\eta - \xi}{|\eta - \xi|} \geq -1/4$,
which implies in particular
\begin{equation*}
\frac{\xi}{|\xi|} \cdot \frac{\eta}{|\eta|} \geq -\frac{1}{4} \, .
\end{equation*}
Therefore, the support of $\chi$ does not include $\mR_{- +}$.  Moreover on  the support of $1 - \chi$, that is
$\frac{\xi}{|\xi|} \cdot \frac{\eta - \xi}{|\eta - \xi|} \leq 1/4$,
frequencies are localized around $\mR_{- +}$  and away from $\mS_{- +} \cap \mT_{- +}^c$.
Let  $
\chi_+ = \chi$ and  $\chi_- = 1 - \chi$
and write 
\begin{equation*}
h_0(h,h) = h_+ (h,h)  +  h_- (h,h)
\end{equation*}
with
\begin{align}
\label{Btheta+}
&\what{h}_\pm (\xi)   := \int_1^t \! \int  e^{is \phi (\xi,\eta)}  \, \nabla_\eta \phi (\xi,\eta) \,
			\chi_\pm (\xi,\eta)
			\what{h} (\eta) \what{h} (\xi-\eta)  \, d\eta \, ds
\end{align}

\subsection{Estimate of $\sup_t t {\left\| e^{i t \Lambda} h_+ (h,h) \right\|}_{L^\infty}$\label{Linftyh+}}
This is the term whose frequencies are away from the time resonant set $\mT$.
Therefore
we can integrate by parts in time to obtain
\begin{subequations}
\begin{align}
\label{Btheta+1}
\what{h}_+ (t,\xi)   =& \left. \int  e^{is \phi}  \frac{\nabla_\eta \phi (\xi,\eta)}{\phi(\xi,\eta)}
			\chi_+ (\xi,\eta)
			\what{h} (\eta) \what{h} (\xi-\eta)  \, d\eta
			\right]_1^t
\\
\label{Btheta+2}
& - \int_1^t \! \int  e^{is \phi}  \frac{\nabla_\eta \phi (\xi,\eta)}{\phi(\xi,\eta)}
			\chi_+ (\xi,\eta)
			\partial_s \left( \what{h} (\eta) \what{h} (\xi-\eta) \right) \, d\eta \, ds \, 
\\
\nn
\overset{def}{=} &\hat b(t) -\hat b(1)  - \int_1^t \! \int  e^{is \phi}  \frac{\nabla_\eta \phi (\xi,\eta)}{\phi(\xi,\eta)}
			\chi_+ (\xi,\eta)
			\partial_s \left( \what{h} (\eta) \what{h} (\xi-\eta) \right) \, d\eta \, ds \, 
\end{align}
\end{subequations}

The bound on  $\sup_t t {\left\| e^{i t \Lambda}b(1) \right\|}_{L^\infty}$ follows from the assumptions on the initial data.  To bound 
  $\sup_t t {\left\| e^{i t \Lambda}b(t) \right\|}_{L^\infty}$, we note that in physical space 
\begin{equation}
\label{eitgradh_+}
e^{i t \Lambda} b(t) =  
		T_{\frac{\nabla_\eta \phi (\xi,\eta)}{\phi(\xi,\eta)} \chi_+ (\xi,\eta) } 
		\left( e^{-it\grad} h , e^{it\grad} h \right) \, .
\end{equation}
By the definition of the cutoff, $\phi$ does not vanish on the support of $\chi_+$.
Nevertheless, dividing by $\phi$ introduces some singularity, which we need to take care of.
Now we claim that 
\begin{equation}
\label{eitgradh_+1}
	e^{i t \Lambda} b(t) \sim T_{ \frac{\mu_0 (\xi,\eta)}{|\eta|} + \frac{\mu_0^\p (\xi,\eta) }{|\xi - \eta|}} 
		\left( e^{-it\grad} h , e^{it\grad} h \right) \, .
\end{equation}
for some $\mu_0, \mu_0^\p \in B_0$, so that this term  is analogous to  $g$  \eqref{ghat},
which was  estimated  in section \ref{Linftyg}.
 To verify the claim first note that
\begin{equation}
\label{1/phi}
\frac{1}{\phi(\xi,\eta)} = -\frac{1}{2}  \frac{|\xi| + |\eta| + |\xi-\eta|}{\xi \cdot \eta + |\xi||\eta|} 
=  \frac{\mu_0 (\xi,\eta)}{|\xi|}  +   \frac{\mu_0^\p (\xi,\eta)}{|\eta|}  \, ,
\end{equation}
where the last identity holds since $\xi \cdot \eta + |\xi| |\eta| \geq \frac{3}{4} |\xi| |\eta|$ 
on the support of $\chi_+$.
Furthermore, since $\nabla_\eta \phi$ vanishes at $\xi = 0$, one can see that
$\frac{\nabla_\eta \phi(\xi,\eta)}{|\xi|} = \frac{\mu_0(\xi,\eta)}{|\eta - \xi|}$.
It then follows
\begin{equation}
\label{dphi/phi}
\frac{\nabla_\eta \phi (\xi,\eta)}{\phi(\xi,\eta)} 
	 = \frac{\mu_0(\xi,\eta) }{|\eta|} +  \frac{\mu_0^\p(\xi,\eta)}{|\xi - \eta|} \, .
\end{equation}

To estimate the remaining term let us denote by $\psi(\xi,\eta) =  \frac{\nabla_\eta \phi(\xi,\eta)}{\phi (\xi,\eta)}  \chi_+ (\xi,\eta)$ and by 
 \begin{equation}
B(t):= e^{it\Lambda}\F^{-1}\int_1^t \! \int  e^{is \phi} \psi (\xi,\eta)
			\partial_s \left( \what{h} (\eta) \what{h} (\xi-\eta) \right) \, d\eta \, ds
\end{equation}
Using the dispersive estimate, and interpolating between weighted $L^2$ norms we see that
\begin{subequations}
\begin{align*}
t{\left\| B(t) \right\| }_{L^\infty} 
		\lesssim 
		&
		\int_1^t \! {\left\|  e^{-is\Lambda}  
		T_{\psi}
		\left( e^{-is\Lambda} \partial_s h,  e^{is\Lambda}  h \right) \right\|}_{\dot{W}^{2,1}}  ds
\lesssim  
		 \int_1^t \! {\left\|  e^{-is\Lambda}  
		T_{|\xi|^2\psi}
		\left( e^{-is\Lambda} \partial_s h,  e^{is\Lambda}  h \right) \right\|}_{1} ds
\\
\lesssim &
\int_1^t \! {\left\| |x| e^{-is\Lambda}  
		T_{|\xi|^2\psi}
		\left( e^{-is\Lambda} \partial_s h,  e^{is\Lambda}  h \right) \right\|}_{2}^{\frac{1}{2}}
{\left\| {|x|}^2 e^{-is\Lambda}  
		T_{|\xi|^2\psi}
		\left( e^{-is\Lambda} \partial_s h,  e^{is\Lambda}  h \right) \right\|}_{2}^{\frac{1}{2}} \, ds \, .
\end{align*}
\end{subequations}
Therefore to obtain the desired $L^\infty$ decay of $t^{-1}$,
it will be sufficient to prove the  following bounds:
\begin{subequations}
\begin{align}
\label{BSx}
& {\left\| x \, e^{-is\Lambda} 
		T_{|\xi|^2\psi}
		\left( e^{-is\Lambda} \partial_s h,  e^{is\Lambda}  h \right) \right\|}_2
			\lesssim  \frac{1}{s^{\frac{7}{4}}} {\| u \|}_X^2
\\
\label{BSx^2}
& {\left\| {|x|}^2 e^{-is\Lambda} 
		T_{|\xi|^2\psi}
		\left( e^{-is\Lambda} \partial_s h,  e^{is\Lambda}  h \right) \right\|}_2
			\lesssim \frac{1}{s^\frac{3}{4}} {\| u \|}_X^2 \, .
\end{align}
\end{subequations}

\vskip15pt
\subsubsection{Proof of \eqref{BSx}.}
We need to look at the different terms  of
\begin{equation}\label{F}
\hat F :=|\xi| \nabla_\xi  \int  e^{is \phi}
		|\xi| \psi(\xi,\eta) 
			\partial_s \what{h} (\eta) \what{h} (\xi-\eta)  \, d\eta \, .
\end{equation}
Using \eqref{xidxiphi} and integrating by parts in frequency, we get
\begin{subequations}
\begin{align}
\label{BSx1}
&  \int s   e^{is \phi}  
		|\xi | \nabla_\eta \phi(\xi,\eta) \chi_+ (\xi,\eta)
			\partial_s \what{h} (\eta) \what{R h} (\xi-\eta)  \, d\eta
\\
\label{BSx2}
& \int  e^{is \phi}  \left[ {|\xi|} \nabla_\xi - |\eta| \nabla_\eta \right] 
		\left(|\xi|\psi(\xi,\eta)  \right)
			\partial_s \what{h} (\eta) \what{h} (\xi-\eta)  \, d\eta
\\
\label{BSx4}
& \int  e^{is \phi} |\xi|\psi(\xi,\eta)
		\, |\eta| \nabla_\eta  \partial_s \what{h} (\eta) \what{h} (\xi-\eta)  \, d\eta
\\
\label{BSx5}
& \int  e^{is \phi} 
		|\xi|\psi(\xi,\eta)
		\partial_s \what{h} (\eta) \, |\eta - \xi| \nabla_\eta  \what{h} (\xi-\eta)  \, d\eta
\end{align}
\end{subequations}
plus similar and easier terms. We now proceed to estimate the above terms.

\subparagraph{Estimate of ${\| \eqref{BSx1} \|}_2$.}
We integrate by parts in frequency using the presence of $\nabla_\eta \phi$.
This gives the following contributions:
\begin{subequations}
\begin{align}
\label{BSx11}
&  \int  e^{is \phi}  |\xi | \mu_0 (\xi,\eta) 
		\partial_s \nabla_\eta \what{h} (\eta) \what{h} (\xi-\eta)  \, d\eta
\\
\label{BSx12}
&  \int  e^{is \phi}  |\xi | \mu_0 (\xi,\eta)
		\partial_s \what{h} (\eta)  \nabla_\eta \what{h} (\xi-\eta)  \, d\eta
\\
\label{BSx13}
&  \int  e^{is \phi}  |\xi | \nabla_\eta \mu_0 (\xi,\eta) 
		\partial_s \what{h} (\eta) \what{h} (\xi-\eta)  \, d\eta
\\
\nn
& + \quad \mbox{``similar terms''}  \, .
\end{align}
\end{subequations}
where, as usual, $\mu_0$ denotes a generic symbol in the class $\B_0$.
\eqref{BSx11} and can be bounded by an
$L^{ \frac{2\d}{2-\d} } \times L^\frac{1}{\d}$ estimate, with $\d$ small enough,
using \eqref{boundxdsf}:
\begin{align*}
{ \left\| \eqref{BSx11} \right\| }_{L^2} & =
	{\left\| T_{|\xi | \mu_0 (\xi,\eta) }
	\left( e^{- i s \Lambda} \partial_s x h ,  e^{i s \Lambda} h \right) \right\| }_{L^2}
\\
& \lesssim 
	{\left\|  e^{- i s \Lambda} \partial_s x h \right\|}_{W^{1, {\frac{2}{1-2\d} }}}
	{\left\|  e^{i s \Lambda} h  \right\|}_{W^{1,\frac{1}{\d}}}
	\\
& \lesssim 
	\frac{1}{s} s^\g s^{\d_{N}} {\| u \|}^2_{X}
	\frac{1}{s^{1-2\d}} s^{\d_{N}} {\| u \|}_{X}
\lesssim
	\frac{1}{s^\frac{7}{4}} {\| u \|}_{X}^3  \, ,
\end{align*}
provided we have chosen $\d$ small enough, and $N$ large enough.
One can estimate \eqref{BSx12} in a similar fashion by
an $L^\frac{1}{\d}  \times L^{ {\frac{2}{1-2\d} }}$,
making use of  \eqref{bounddsf}.

Since $\nabla_\eta \mu_0 =  \frac{\mu_0^{\p}}{|\eta|} + \frac{\mu_0^{\p\p}}{|\xi-\eta|}$,
for some $\mu_0^{\p}, \mu_0^{\p\p} \in \B_0$, then \eqref{BSx13} can be written  as 

\begin{equation}
\label{BSx131}
 e^{is \Lambda}  T_{ |\xi| \mu_0^\p (\xi,\eta) }
		\left( \frac{1}{\grad} e^{-is\Lambda} \partial_s h,  e^{is\Lambda} h \right)
+  e^{is \Lambda}  T_{  |\xi| \mu_0^{\p\p} (\xi,\eta)   }
		\left( e^{-is\Lambda} \partial_s h, \frac{1}{\Lambda} e^{is\Lambda} h  \right) = I + II  \, .
\end{equation}
By Theorem \eqref{theoop} and Sobolev's embeddings we can estimate
\begin{align*}
   {\left\| I \right\| }_{L^2}
& \lesssim
	{\left\| \frac{1}{\Lambda} e^{- i s \Lambda} \partial_s h\right\|}_{W^{1,6}}
	{\left\|  e^{is\Lambda} h  \right\|}_{W^{1,3}}
\lesssim {\left\|  e^{i s \Lambda} \partial_s h \right\|}_{H^1}
	\frac{1}{s^\frac{1}{3}} s^{\d_{N}} {\| u \|}_X
\\
& \lesssim \frac{1}{s^{2-}} s^\g s^{\d_{N}}  {\| u \|}_X^2  
	\frac{1}{s^\frac{1}{3}} s^{\d_{N}} {\| u \|}_X
\lesssim \frac{1}{s^\frac{7}{4}} {\| u \|}_X^3 \, ,
\end{align*}
if $N$ is large enough so that $ 2\d_{N} < \frac{1}{4}$.  $II$ can be estimated 
in a similar fashion.

Next we estimate  \eqref{BSx2}.   
Using successively \eqref{xidxiphi} and \eqref{dphi/phi}, one can see that
\begin{align*}
\left[ {|\xi|} \nabla_\xi - |\eta| \nabla_\eta \right]  
	\left( |\xi| \psi  (\xi,\eta) \right)
	& =  \mu_0 (\xi,\eta) \frac{|\xi| \nabla_\eta \phi(\xi,\eta)}{\phi(\xi,\eta)}
 = |\xi| \frac{\mu_0^\p (\xi,\eta) }{|\eta|} +  |\xi| \frac{\mu_0^{\p\p} (\xi,\eta)}{|\xi - \eta|} \, .
\end{align*}
Therefore, the contribution coming from \eqref{BSx2} is identical to
those of \eqref{BSx131} that we have just estimated above.

To estimate  \eqref{BSx4}  we note that 
\begin{align*}
\F^{-1} \eqref{BSx4} & =  e^{- is \Lambda}
	T_{ \mu_1 (\xi,\eta) }
	\left(  e^{-is \Lambda}  \partial_s x h , e^{is \Lambda} h \right)
\end{align*}
since $\frac{|\xi||\eta|}{\phi(\xi,\eta)} \in \B_1$   by  \eqref{1/phi}.  This 
 can be estimated using lemma \ref{lem6}
\begin{align*}
{ \left\| \eqref{BSx4} \right\| }_{L^2} & \lesssim 
		{\left\|  e^{- i s \Lambda} \partial_s x h \right\|}_{W^{1, {\frac{2}{1-2\d} }}}
		{\left\|  e^{i s \Lambda}  h  \right\|}_{W^{1,\frac{1}{\d}}}
 \lesssim	{\left\|  \partial_s  x h \right\|}_{H^2}
		\, \frac{1}{s} \, s^{2 \d  + \d_{N}} {\left\|  u \right\|}_{X} 
\\
& \lesssim \frac{1}{s} s^{\g + \d_{N}} {\| u \|}_X^2
		\, \frac{1}{s} \, s^{2 \d  + \d_{N}} {\left\|  u \right\|}_{X} 
\lesssim \frac{1}{s^\frac{7}{4}} {\| u \|}_X^3 \, ,
\end{align*}
provided $N$ is large enough.

Finally,  \eqref{BSx5}  is given by
\begin{equation}
\label{BSx51}
 e^{ -is \Lambda  }  T_{ |\xi| \mu_0 (\xi,\eta)}
	\left( \frac{1}{\Lambda} e^{-is\Lambda} \partial_s h , e^{is\Lambda} x \Lambda h \right)
 +  e^{- is \Lambda  }  T_{ |\xi| \mu_0 (\xi,\eta)}
	\left(  e^{-is\Lambda} \partial_s h , e^{is\Lambda} x  h \right)
\end{equation}
These terms are similar to each other and can be  bounded  using  Theorem \ref{theoop} and estimate \eqref{bounddsf}:
\begin{align*}
   {\left\| T_{\mu_0 (\xi,\eta)  }
	\left( \frac{e^{ -is\Lambda} }{\Lambda} \partial_s h , e^{is\Lambda} x \Lambda h  \right) \right\| }_{\dot{H}^1}\!\!
 \lesssim &
		{\left\|  e^{- is \Lambda}  \partial_s h  \right\|}_{L^6}
		{\left\|  e^{is\Lambda} x \Lambda h  \right\|}_{L^3}
+
		{\left\|  \frac{e^{- is \Lambda} }{\Lambda}  \partial_s h  \right\|}_{W^{1,6}}
		{\left\|  e^{is\Lambda} x {\Lambda}^2 h  \right\|}_{L^2}
\\
 \lesssim &
		{\left\|  \partial_s h \right\|}_{H^1} {\left\|  x \Lambda h  \right\|}_{H^1}
\lesssim \frac{1}{s^{2-}} s^{\g + \d_{N}} {\| u \|}_X^3
\lesssim \frac{1}{s^\frac{7}{4}} {\| u \|}_X^3 \, ,
\end{align*}
provided, as usual, $N$ is chosen large enough.

\subsubsection{Proof of \eqref{BSx^2}.}   
In order to prove \eqref{BSx^2} we need to estimate the $L^2$ norm of 
$\nabla_\xi\hat  F$, where $\hat F$ is given by \eqref{F}.   First we integrate by parts in frequency  to reduce the bilinear term to prototypes given by  \eqref{BSx1}-\eqref{BSx5}, then we apply $\nabla_\xi$.  Here we note that  if  $\nabla_\xi$ hits  the phase   we lose   a factor of  $s$.
Therefore, each time we differentiate $e^{is\phi}$ in \eqref{BSx1}-\eqref{BSx5}, we can perform the same estimates  done in the previous section and get the stated  bound $s \, s^{-\frac{7}{4}}$.   Thus the only care we need to take is when the $\nabla_\xi$ hits the bilinear symbol,
or where $\nabla_\xi$ hits the profile $\what{h}(\xi-\eta)$.   These are the terms we will explicitly show how to bound.

From the discussion above we can reduce   ${\| \nabla_\xi \eqref{BSx1} \|}_2$ to
\begin{subequations}
\begin{align}
\label{BSxx12}
&  \int s e^{is \phi}  |\xi | \mu_0 (\xi,\eta)
		\partial_s \what{h} (\eta)  \frac{ \what{h} (\xi-\eta)}{|\xi-\eta|}  \, d\eta
\\
\label{BSxx13}
&  \int s e^{is \phi}  |\xi |  \mu_0 (\xi,\eta) 
		\partial_s \what{h} (\eta)  \nabla_\xi \what{h} (\xi-\eta)  \, d\eta
\end{align}
\end{subequations}
where, as usual, $\mu_0$ denotes a generic symbol in the class $\B_0$.
%
%
%
%
Now, observe that \eqref{BSxx12} and  \eqref{BSxx13},
are  analogous to \eqref{BSx13} and \eqref{BSx12},respectively, with an extra factor of $s$.
This shows immediately that \eqref{BSxx12} and \eqref{BSxx13} satisfy $L^2$ bounds of order $s^{-\frac{3}{4}}$,
whence $\nabla_\xi \eqref{BSx1}$ does too.

Let us first recall that \eqref{BSx2} is given by the following two terms:
\[
 \int  e^{is \phi}  |\xi| \mu_0(\xi,\eta)
			\frac{ \partial_s \what{h} (\eta)}{|\eta|}  \what{h} (\xi-\eta)  \, d\eta,
\qquad  \int  e^{is \phi}  |\xi| \mu_0(\xi,\eta)
			\partial_s \what{h} (\eta) \frac{\what{h} (\xi-\eta)}{|\xi-\eta|}  \, d\eta
\]
Applying $\nabla_\xi$ to these terms, and disregarding the terms when $\nabla_\xi$ hits  the oscillating phase,  we get
\begin{subequations}
\begin{align}
\label{BSxx21}
& \int  e^{is \phi}  |\xi| \mu_0(\xi,\eta)
			\frac{\partial_s \what{h} (\eta)}{|\eta|}  \nabla_\xi \what{h} (\xi-\eta)  \, d\eta
\\
\label{BSxx22}
& \int  e^{is \phi}  |\xi| \mu_0(\xi,\eta)
			\partial_s \what{h} (\eta) \frac{\what{h} (\xi-\eta)}{{|\xi-\eta|}^2}  \, d\eta
\end{align}
\end{subequations}
plus other ``similar or easier'' contributions.
The first term above is very similar to \eqref{BSx51} and can therefore be bounded in an analogous fashion.
The second term can be bounded by means of Theorem \ref{theoop}, Sobolev's embedding,
lemma \ref{lem6}   as follows:
\begin{align*}
{ \left\| \eqref{BSxx22} \right\| }_{L^2} & =
	{\left\| T_{|\xi | \mu_0 (\xi,\eta) }
	\left( e^{- i s \Lambda} \partial_s h,  \frac{ e^{i s \Lambda}}{\grad^2} h \right) \right\| }_{L^2}
 \lesssim 
	{\left\|  e^{- i s \Lambda} \partial_s h \right\|}_{W^{1,3}}
	{\left\|   \frac{e^{i s \Lambda}}{\grad^2} h  \right\|}_{W^{1,6}}
\\
& \lesssim 
	{\left\|  \partial_s h \right\|}_{H^2}
	{\left\|  \langle x \rangle h  \right\|}_{L^2}
 \lesssim 
	\frac{1}{s} s^{\d_{N} + \g} {\| u \|}_{X}^3 
\lesssim  \frac{1}{s^\frac{3}{4}}  {\| u \|}_{X}^3 \, .
\end{align*}

Next we consider
\begin{equation*}
\eqref{BSx4}  =  \int e^{is \phi} \mu_1 (\xi,\eta)
			\partial_s \nabla_\eta \what{h} (\eta) \what{h} (\xi-\eta) \, d\eta
\end{equation*}
and apply  $\nabla_\xi$ to it.   We get only one contribution which differs from the ones 
which have been previously estimated, namely
\begin{equation}
\label{BSxx41}
\int e^{is \phi} \mu_1 (\xi,\eta)
			\partial_s \nabla_\eta \what{h} (\eta) \nabla_\xi \what{h} (\xi-\eta) \, d\eta \, .
\end{equation}
This term can be treated again using Theorem \eqref{theoop}, Sobolev's embedding,
and lemma \ref{lem6}:
\begin{align*}
{ \left\| \eqref{BSxx41} \right\| }_{L^2} & =
	{\left\| T_{ \mu_1 (\xi,\eta) }
	\left( e^{- i s \Lambda} \partial_s x h,  e^{i s \Lambda} x h \right) \right\| }_{L^2}
 \lesssim 
	{\left\|  e^{- i s \Lambda} \partial_s x h \right\|}_{W^{1,6}}
	{\left\|  e^{i s \Lambda} x h  \right\|}_{W^{1,3}}
\\
& \lesssim 
	{\left\|  \partial_s x h \right\|}_{H^2}
	{\left\|  x h  \right\|}_{H^2}
\lesssim 
	\frac{1}{s} s^{\d_{N} + \g}  {\| u \|}_{X}^3 
\lesssim \frac{1}{s^\frac{3}{4}}  {\| u \|}_{X}^3  \, ,
\end{align*}
provided $N$ is large enough.

Finally  applying $\nabla_\xi$ to \eqref{BSx5}, there are only
two contributions that differ from the ones previously obtained and already estimated, which are
\begin{subequations}
\begin{align}
\label{BSxx51}
& \int e^{is \phi} \mu_1 (\xi,\eta)
			\frac{ \partial_s \what{h} (\eta)}{|\eta|} 
			|\xi-\eta| \nabla_\eta \nabla_\xi \what{h} (\xi-\eta) \, d\eta
\\
\label{BSxx52}
& \int e^{is \phi} \mu_1 (\xi,\eta)
			\partial_s \what{h} (\eta) 
			\nabla_\eta \nabla_\xi \what{h} (\xi-\eta) \, d\eta \, .
\end{align}
\end{subequations}
As in the previous estimates, we can bound
\begin{align*}
{ \left\| \eqref{BSxx51} \right\| }_{L^2} & =
	{\left\| T_{ \mu_1 (\xi,\eta) }
	\left( e^{- is\Lambda} \frac{1}{\grad} \partial_s h, e^{is\Lambda} \grad {|x|}^2 h \right) \right\| }_{L^2}
\\
& \lesssim 
	{\left\|  e^{- is\Lambda} \frac{1}{\grad} \partial_s h \right\|}_{W^{1,6}}
	{\left\|  e^{is\Lambda} \grad {|x|}^2 h \right\|}_{H^1}
 +
	{\left\|  e^{- is\Lambda} \partial_s h \right\|}_{L^6}
	{\left\|  e^{is\Lambda} \grad {|x|}^2 h \right\|}_{L^3}
\\
& \lesssim 
	{\left\|  \partial_s h \right\|}_{H^1}
	{\left\|  \grad {|x|}^2 h  \right\|}_{H^1}
 \lesssim 
	\frac{1}{s} {\| u \|}_{X}^2  s^a {\| u \|}_{X} 
\lesssim \frac{1}{s^\frac{3}{4}}  {\| u \|}_{X}^3  \, .
\end{align*}
Similarly, using \eqref{fractional}, we also have
\begin{align*}
{ \left\| \eqref{BSxx52} \right\| }_{L^2} & =
	{\left\| T_{ \mu_1 (\xi,\eta) }
	\left( e^{- is\Lambda} \partial_s h, e^{is\Lambda}  {|x|}^2 h \right) \right\| }_{L^2}
 \lesssim 
	{\left\| e^{-is\Lambda} \partial_s h \right\|}_{W^{1,3}}
	{\left\|  e^{is\Lambda} {|x|}^2 h \right\|}_{W^{1,6}}
\\
& \lesssim 
	{\left\|  \partial_s h \right\|}_{H^2}
	{\left\|  \grad {|x|}^2 h  \right\|}_{H^1}
\lesssim 
	\frac{1}{s} {\| u \|}_{X}^2  s^a {\| u \|}_{X} 
\lesssim \frac{1}{s^\frac{3}{4}}  {\| u \|}_{X}^3  \,  .
\end{align*}
This concludes the proof of \eqref{BSx^2}, and hence shows the desired $L^\infty$ bound on $h_+$.

\subsection{Estimate of $\sup_t t {\| e^{i t \Lambda} h_- (h,h) \|}_{L^\infty} $\label{Linftyh-}}
This is the term whose frequencies are localized around the resonant set $\mR$, 
and, therefore, is the hardest to treat.
By definition, frequencies in the support of $1 - \chi_+$ are away from $\mS \cap \mT^c$
and, in particular, satisfy
\begin{equation}
\xi \cdot (\xi - \eta) + |\xi| |\xi - \eta| \geq \frac{3}{4} |\xi| |\xi - \eta|
\qquad
\mbox{or, equivalently}
\qquad
\cos(\xi, \xi - \eta) \geq - \frac{1}{4} \, .
\end{equation}
This fact allows us to exploit the space-resonance of the phase through the following identity:
\begin{equation}
\label{phigradphi}
\phi = \frac{(|\xi| + |\eta| + |\xi-\eta|) (|\eta||\xi-\eta|) }{\xi \cdot (\xi - \eta) + |\xi||\xi-\eta|} 
							{|\nabla_\eta \phi|}^2  \, .
\end{equation}
Thanks to the above identity and \eqref{xidxiphi}, we can express $\nabla_\xi \phi$
in terms of $\nabla_\eta \phi$.
This will introduce some singularity which needs to be carefully analyzed,
but, eventually, it will imply good weighted $L^2$-estimates for $h_-$.
From these we will deduce the $L^\infty$ decay of $e^{it\Lambda} h_-$,
similarly to how we did for $h_+$.

In real space
\begin{equation*}
e^{it\Lambda} h_-(h,h)   =  e^{it\Lambda} \int_1^t  e^{-is\Lambda} 
	T_{ \nabla_\eta \phi (\xi,\eta) \chi_- (\xi,\eta) } \left( e^{-is\grad} h , e^{is\grad} h \right) \, ds \, ,
\end{equation*}
whence, using the dispersive estimate,
\begin{subequations}
\begin{align*}
{\left\| e^{it\Lambda} h_- \right\|}_\infty &  \lesssim
		\frac{1}{t} \int_1^t {\left\| e^{-is\Lambda} 
		T_{ \nabla_\eta \phi (\xi,\eta) \chi_- (\xi,\eta) } \left( e^{-is\grad} h , e^{is\grad} h \right)
			 \right\|}_{\dot{W}^{2,1}}  \, ds
\\
& \lesssim 
\frac{1}{t} \int_1^t {\left\| e^{-is\Lambda} 
		T_{ {|\xi|}^2 \nabla_\eta \phi (\xi,\eta) \chi_- (\xi,\eta) } \left( e^{-is\grad} h , e^{is\grad} h \right) 
			\right\|}_1  \, ds
\\
& \lesssim 
\frac{1}{t} \int_1^t {\left\| {|x|}^2 e^{-is\Lambda} 
		T_{ {|\xi|}^2 \nabla_\eta \phi (\xi,\eta) \chi_- (\xi,\eta) } \left( e^{-is\grad} h , e^{is\grad} h \right)
			 \right\|}_2^{\frac{1}{2}}
\\
& \quad \times {\left\| x \, e^{-is\Lambda} 
		T_{ {|\xi|}^2 \nabla_\eta \phi (\xi,\eta) \chi_- (\xi,\eta) } \left( e^{-is\grad} h , e^{is\grad} h \right)
			 \right\|}_2^{\frac{1}{2}} \, ds
\end{align*}
\end{subequations}
Therefore, in order to prove ${\| e^{it\Lambda} h_- \|}_\infty \lesssim \frac{1}{t}$
it will suffice to establish the two following estimates:
\begin{subequations}
\begin{align}
\label{BRx}
& {\left\| x \, e^{-is\Lambda} 
	T_{ {|\xi|}^2 \nabla_\eta \phi (\xi,\eta) \chi_- (\xi,\eta) } \left( e^{-is\grad} h , e^{is\grad} h \right) \right\|}_2
			\lesssim \frac{1}{s^\frac{7}{4}} {\| u \|}_X^2
\\
\label{BRx^2}
& {\left\| {|x|}^2 e^{-is\Lambda} 
	T_{ {|\xi|}^2 \nabla_\eta \phi (\xi,\eta) \chi_- (\xi,\eta) } \left( e^{-is\grad} h , e^{is\grad} h \right) \right\|}_2
			\lesssim \frac{1}{s^\frac{3}{4}} {\| u \|}_X^2 \,  .
\end{align}
\end{subequations}

As usual, we try to bound in $L^2$ the expression
\begin{equation*}
{|\xi|}^2 \nabla_\xi  \int  e^{is \phi}  \nabla_\eta \phi (\xi,\eta) \chi_- (\xi,\eta) 
		\what{h} (\eta) \what{h} (\xi-\eta)  \, d\eta  \, .
\end{equation*}
Using \eqref{xidxiphi} as already done previously,  this produces the following terms
\begin{subequations}
\begin{align}
\label{BRx1}
& \int s  e^{is \phi} |\xi|  \phi (\xi,\eta)
		\nabla_\eta \phi (\xi,\eta) \chi_- (\xi,\eta) 
		\what{h} (\eta) \what{Rh} (\xi-\eta)  \, d\eta
\\
\label{BRx2}
& \int  e^{is \phi} |\xi| \left[ |\xi| \nabla_\xi - |\eta| \nabla_\eta  \right] 
		\left[ \nabla_\eta \phi (\xi,\eta) \chi_- (\xi,\eta) \right] 
		\what{h} (\eta) \what{h} (\xi-\eta)  \, d\eta
\\
\label{BRx3}
& \int  e^{is \phi} |\xi| \nabla_\eta \phi (\xi,\eta) \chi_- (\xi,\eta) 
		\what{h} (\eta) \left[ |\xi| \nabla_\xi - |\eta| \nabla_\eta  \right] \what{h} (\xi-\eta)  \, d\eta
\\
\label{BRx4}
& \int  e^{is \phi}  |\xi|  \nabla_\eta \phi (\xi,\eta) \chi_- (\xi,\eta) 
		|\eta| \nabla_\eta  \what{h} (\eta)  \what{h} (\xi-\eta)  \, d\eta
\end{align}
\end{subequations}
plus ``similar and easier terms''.  Notice that the first term contains a factor of $s$, and has  two symbols vanishing:
$\phi$ and  $\nabla_\eta \phi$. Also, all of the terms \eqref{BRx2}-\eqref{BRx4} contain a$\nabla_\eta \phi$ symbol.

To estimate  \eqref{BRx1}, we have from \eqref{phigradphi} 
\begin{equation}
\label{phigradphi2}
{|\xi|} \phi = \frac{ (|\xi| + |\eta| + |\xi-\eta|) |\eta|}{1 + \cos(\xi, \xi - \eta)} 
		{|\nabla_\eta \phi|}^2
		= \mu_1 (\xi,\eta)  |\eta| {|\nabla_\eta \phi|}^2   \, ,
\end{equation}
so that
\begin{equation}
\label{BRx10}
\eqref{BRx1} = \int s \, e^{is \phi} \mu_1 (\xi,\eta)  |\eta| {|\nabla_\eta \phi(\xi,\eta)|}^2 
		\nabla_\eta \phi (\xi,\eta) \chi_- (\xi,\eta) 
		\what{h} (\eta) \widehat{Rh} (\xi-\eta)  \, d\eta \, .
\end{equation}
Before integrating by parts in $\eta$,
notice that on the support of $\chi_-$, we must have $|\eta| , |\xi| \leq 2 |\eta - \xi|$.
This can be easily seen from the following:
\begin{equation*}
{|\xi-\eta|}^2 = {|\xi|}^2 + {|\eta|}^2 - 2 \xi \cdot \eta \geq 
{|\xi|}^2 + {|\eta|}^2 - \frac{1}{2} |\xi| |\eta|  \geq \frac{3}{2} |\xi| |\eta| \, ,
\end{equation*}
so that $\min \{ |\xi|, |\eta| \} \leq |\xi-\eta|$, whence $|\xi|, |\eta| \leq 2 |\eta - \xi|$.
Thus, on the support of $\chi_-$, symbols in $\B_1$ behave essentially like $|\xi-\eta|$.
Now we integrate by parts in \eqref{BRx10} obtaining a term
\begin{subequations}
\begin{align}
\label{BRx11}
& \int  e^{is \phi} \mu_1 (\xi,\eta) \chi_- (\xi,\eta) \nabla_\eta \phi (\xi,\eta)
	 	|\eta| \nabla_\eta \what{h} (\eta) \what{h} (\xi-\eta)  \, d\eta
\end{align}
\end{subequations}
plus ``similar or easier'' ones.
Notice that even after this integration by parts,
a symbol like $\nabla_\eta \phi$ survives, so that we can integrate by parts once again.
This gives the following contributions:
\begin{subequations}
\begin{align}
\label{BRx14}
& \frac{1}{s} \int  e^{is \phi} \mu_1 (\xi,\eta) \chi_- (\xi,\eta)
		\, |\eta| \nabla_\eta^2 \what{h} (\eta) \what{h} (\xi-\eta)  \, d\eta
\\
\label{BRx15}
& \frac{1}{s} \int  e^{is \phi} \mu_1 (\xi,\eta) \chi_- (\xi,\eta)
	\, |\eta| \nabla_\eta \what{h} (\eta) \nabla_\eta \what{h} (\xi-\eta)  \, d\eta \, ,
\end{align}
\end{subequations}
plus other ``symmetric and easier terms''.
As already done before, we let $0 < \d \ll 1$ and invoke Theorem \ref{theoop} to obtain
\begin{align*}
{ \left\| \eqref{BRx14} \right\| }_{L^2} & = \frac{1}{s}
	{\left\| T_{ \mu_1 (\xi,\eta) }
	\left( e^{- is\Lambda} \grad {|x|}^2 h, e^{is\Lambda} h \right) \right\| }_{L^2}
 \lesssim \frac{1}{s}
	{\left\| e^{-is\Lambda} \grad {|x|}^2 h \right\|}_{H^1}
	{\left\|  e^{is\Lambda} h \right\|}_{W^{2,\frac{1}{\d}}}
\\
& \lesssim \frac{1}{s}
	s^a {\| u \|}_{X}
	\frac{1}{s^{1-2\d}} \, s^{\d_{N}} {\| u \|}_{X} 
\lesssim 
	\frac{1}{s^\frac{7}{4}}  {\| u \|}_{X}^2  \,  ,
\end{align*}
provided $\d$ is chosen small enough, and $N$ large enough.
Using $\mu_1 (\xi,\eta) = \mu_0(\xi,\eta) {|\xi-\eta|}$ on the support of $\chi_-$,
and, again, Theorem \ref{theoop} and dispersive estimates, we have
\begin{align*}
{ \left\| \eqref{BRx15} \right\| }_{L^2} & = \frac{1}{s}
	{\left\| T_{ \mu_1 (\xi,\eta) }
	\left( e^{- is\Lambda} \grad {|x|} h, e^{is\Lambda} {|x|} h \right) \right\| }_{L^2}
 \lesssim \frac{1}{s}
	{\left\| e^{-is\Lambda} \grad {|x|} h \right\|}_{L^4}^2
\\
& \lesssim 
	\frac{1}{s^2} {\left\| \grad^2 {\langle x \rangle}^2 h \right\|}_{L^2}^2
 \lesssim \frac{1}{s^2} s^{2b} {\| u \|}_{X}^2
	\lesssim
\frac{1}{s^\frac{7}{4}}  {\| u \|}_{X}^2  \,  .
\end{align*}

In order to estimate \eqref{BRx2} we use subsequently \eqref{xidxiphi} and \eqref{phigradphi2} to obtain
\begin{align*}
 |\xi| \left[ |\xi| \nabla_\xi - |\eta| \nabla_\eta  \right] 
		\left[ \nabla_\eta \phi (\xi,\eta) \chi_- (\xi,\eta) \right]
& = \mu_1(\xi,\eta) \nabla_\eta \phi (\xi,\eta) + \mu_0(\xi,\eta) \frac{|\xi| \phi(\xi,\eta)}{|\xi - \eta|}
\\
& = \mu_1(\xi,\eta) \nabla_\eta \phi (\xi,\eta) \, .
\end{align*}
Therefore, \eqref{BRx2} is of the form
\begin{equation}
\label{BRx21}
\int  e^{is \phi} \mu_1 (\xi,\eta) \nabla_\eta \phi(\xi,\eta)
		\what{h} (\eta) \what{h} (\xi-\eta)  \, d\eta \, ,
\end{equation}
and we can integrate by parts in $\eta$ gaining decay in $s$.
This will give us terms which are easier to estimate than \eqref{BRx14}.
We skip the details.

Since on the support of $\chi_-$, symbols in $\B_1$ behave  like $|\xi-\eta|$,
we can write \eqref{BRx3} as
\begin{align}
\label{BRx31}
& \int  e^{is \phi}  \mu_1 (\xi,\eta) \nabla_\eta \phi(\xi,\eta)
		\what{h} (\eta) |\xi-\eta| \nabla_\eta \what{h} (\xi-\eta)  \, d\eta \, .
\end{align}
This term is analogous to \eqref{BRx11} and can therefore be treated in exactly the same way.
Finally, notice that also \eqref{BRx4}
is of the form \eqref{BRx11} and can be treated in the same way.

\subsubsection{Proof of \eqref{BRx^2}}
In the previous paragraphs we saw that the contribution of the terms \eqref{BRx1}-\eqref{BRx4}
essentially reduces to a term like \eqref{BRx11}, that is
\begin{align}
\label{BRxx0}
& \int  e^{is \phi}  \mu_1 (\xi,\eta) \nabla_\eta \phi(\xi,\eta)
		|\eta| \nabla_\eta \what{h} (\eta)  \what{h} (\xi-\eta)  \, d\eta \, .
\end{align}
In order to establish \eqref{BRx^2} we need to estimate the $L^2$ norm of $\nabla_\xi$
applied to the above term.
Applying $\nabla_\xi$ to \eqref{BRxx0}, will give three terms:

\paragraph{$\nabla_\xi$ hits the phase:}
This will give a term analogous to \eqref{BRxx0} with an extra factor of $s$ in front.
Since  ${\| \eqref{BRxx0} \|}_{L^2} \lesssim s^{-\frac{7}{4}}$, it immediately follows
that this term satisfies a bound of order $s^{-\frac{3}{4}}$.

\paragraph{$\nabla_\xi$ hits the symbol:}
This gives a term of the form
\begin{align*}
& \int  e^{is \phi}  \mu_0 (\xi,\eta)
		|\eta| \nabla_\eta \what{h} (\eta)  \what{h} (\xi-\eta)  \, d\eta \, .
\end{align*}
which can again be treated by Theorem \eqref{theoop} with an a $L^2 \times W^{1,\frac{1}{\d}}$ estimate.

\paragraph{$\nabla_\xi$ hits the profile $\what{f}(\xi-\eta)$:}
This produces a term analogous to \eqref{BRx15} multiplied by a factor of $s$.
A bound of order $s^{-\frac{3}{4}}$ then follows from the estimates already performed.
This concludes the proof of \eqref{BRx^2} yielding the desired $L^\infty$ estimate for $h_- (h,h)$.

\section{Weighted estimates on $h_0 (g,f)$ and $h_1$\label{weightedh_0fg}}

\subsection{Reduction to cubic terms}
From the definition of $h_0$ in \eqref{h0hat} and $g$ in \eqref{ghat} we have
\begin{align}
\nn
\what{h}_0 (g,f) & = \int_1^t \! \int e^{is\phi(\xi,\eta)} \nabla_\eta \phi(\xi,\eta)
		\what{g}(s,\eta) \what{f}(s,\xi-\eta) \, d\eta ds
\\
\label{h0gf}
& = \int_1^t \! \iint e^{is\phi(\xi,\eta,\s)} \nabla_\eta \phi(\xi,\eta) 
		\frac{\mu_0(\eta,\s)}{|\s|}
		\what{f}(s,\s) \what{f}(s,\eta-\s) \what{f}(s,\xi-\eta) \, d\eta d\s ds
\end{align}
where the phase $\phi(\xi,\eta,\s)$ can be any of the combinations
\begin{equation}
\label{cubicphi}
\phi_{+,\pm,\pm} (\xi,\eta,\s) = - |\xi| + |\xi-\eta| \pm |\s| \pm |\eta-\s| \, ,
\end{equation}
and we omitted complex conjugates on the profiles as they play no role in our analysis.
From the definition of $h_1$ in \eqref{h1hat} we have
\begin{align}
\nn
\what{h}_1 (t,\xi) & = \int_1^t \! \int e^{is\phi(\xi,\eta)} \frac{1}{|\eta|}
		\partial_s \left( \what{f}(s,\eta) \what{f}(s,\xi-\eta)  \right) \, d\eta ds
\\
\label{h1hat1}
& = \int_1^t \! \iint e^{is\phi(\xi,\eta,\s)}  \frac{\mu_0(\eta,\s)}{|\eta|}
		\what{f}(t,\s) \what{f}(t,\eta-\s) \what{f}(t,\xi-\eta) \, d\eta d\s ds 
\\
\label{h1hat2}
& + \int_1^t \! \iint e^{is\phi(\xi,\eta,\s)}  \frac{\mu_0(\eta,\s) }{|\xi-\eta|}
		\what{f}(s,\xi-\eta) \what{f}(s,\s-\eta) \what{f}(s,\s) \, d\eta d\s ds 
\end{align}
where $\phi(\xi,\eta,\s)$ is the same as in \eqref{cubicphi}.
Up to relabeling variables, \eqref{h1hat2} is equivalent to \eqref{h0gf}.
Both of these terms are easier to treat than \eqref{h1hat1},
since the singularity in the symbol just pairs with one of the profiles,
while this is not the case for \eqref{h1hat1}.
Therefore, in what follows, we consider that 
\begin{equation*}
\what{h}_0(g,f) + \what{h}_1(f,f) \sim \what{H} (f,f,f) 
\end{equation*}
with $\what{H} = \eqref{h1hat1}$, or, in real space,
\begin{equation}
\label{H}
H (f,f,f) := \int_1^t e^{-is\grad} 
			T_{\frac{\mu_0 (\eta,\s)}{|\eta|}} \left( e^{\pm is \grad} f, e^{\pm is \grad} f \right) 
			\, e^{is \grad} f \, ds \, .
\end{equation}
As already done for $h_0$, we can consider $H$ as a trilinear function of $f = f_0 + g +h$.
Then, we can distinguish between two main contributions: one where all arguments are given by $h$,
and another one where at least one of the arguments is $g$.
In other words, we can regard $H(f,f,f)$ as being given by $H(g,f,f) + H(h,h,h)$.
Since $g$ is quadratic in $f$, $H(g,f,f)$ will be a quartic term, 
with two singularities of the type ${\grad}^{-1}$.
We first show how to treat $H(h,h,h)$ and then indicated how to deal with $H(g,f,f)$.

\subsection{Weighted estimates on $H(h,h,h)$}
In what follows we are going to perform weighted estimates on $H(h,h,h)$.
The $L^\infty$ bound on $e^{-it\grad} H(h,h,h)$ will follow directly from these.
It is important to notice that a formula similar to \eqref{xidxiphi} holds
for the cubic phases\footnote{
As for quadratic phases, one can check that $\nabla_\xi \phi$ vanishes on the space-time resonant 
set. Therefore, an identity like \eqref{xidxiphi3} should not be completely surprising.
Nevertheless, it is not always the case that such a formula will hold in general.
In the case of the wave equation, \eqref{xidxiphi3} is related to scaling invariance.
} \eqref{cubicphi}.
In particular, if we let $\e_1,\e_2 = \pm 1$ we have
\begin{align}
\nn
|\xi| \nabla_\xi \phi_{+,\e_1,\e_2} (\xi,\eta,\s) & = \frac{\eta-\xi}{|\eta-\xi|} \phi_{+,\e_1,\e_2} (\xi,\eta,\s) 
	- (\e_1 |\s| + \e_2 |\eta-\s| ) \nabla_\eta \phi_{+,\e_1,\e_2} (\xi,\eta,\s)
\\
\label{xidxiphi3}
& - \e_1 |\s| \nabla_\s \phi_{+,\e_1,\e_2} (\xi,\eta,\s) \, .
\end{align}

In Fourier space $H(h,h,h)$ is given by
\begin{equation}
\label{Hhat}
\what{H}(h,h,h) = \int_1^t \! \iint e^{is\phi(\xi,\eta,\s)}  \frac{\mu_0(\eta,\s)}{|\eta|}
		\what{h}(t,\s) \what{h}(t,\eta-\s) \what{h}(t,\xi-\eta) \, d\eta d\s ds 
\end{equation}
where $\phi$ can be any of the phases \eqref{cubicphi}.
For simplicity we assume $\phi = \phi_{+,+,+}$ as this will not have any impact on our computations.
We skip the estimate of ${\| H (h,h,h) \|}_{L^2}$ as this easy to show,
and move on to estimate weighted norms of $H(h,h,h)$.

\subsubsection{Estimate of $\sup_t t^{-\g} {\| x H (h,h,h) \|}_{L^2}$}
Applying $\nabla_\xi$ to \eqref{Hhat} gives the following two contributions:
\begin{subequations}
\begin{align}
\label{xHhat1}
& \int_1^t \! \iint s \, e^{is\phi(\xi,\eta,\s)}  \frac{\nabla_\xi\phi(\xi,\eta,\s) \mu_0(\eta,\s)}{|\eta|}
		\what{h}(t,\s) \what{h}(t,\eta-\s) \what{h}(t,\xi-\eta) \, d\eta d\s ds 
\\
\label{xHhat2}
& \int_1^t \! \iint e^{is\phi(\xi,\eta,\s)}  \frac{\mu_0(\eta,\s)}{|\eta|}
		\what{h}(t,\s) \what{h}(t,\eta-\s) \nabla_\xi \what{h}(t,\xi-\eta) \, d\eta d\s ds \, . 
\end{align}
\end{subequations}
Now observe that
\begin{align}
\label{dxiphi/eta}
\frac{\nabla_\xi \phi(\xi,\eta,\s)}{|\eta|}  
	& = \left( - \frac{\xi}{|\xi|} + \frac{\xi -\eta}{|\xi-\eta|} \right) \frac{1}{|\eta|}
	= \frac{\mu_0 (\xi,\eta)}{|\xi-\eta|} \, ,
\end{align}
so that
\begin{align*}
{ \left\| \eqref{xHhat1} \right\| }_{L^2} & = 
	{\left\|  \int_1^t s \, e^{-is\Lambda} T_{ \mu_0 (\xi,\eta) } \left(
	T_{ \mu_0 (\eta,\s) }
	\left( e^{is\Lambda} h, e^{is\Lambda} h \right) ,  e^{is\Lambda} \frac{h}{\grad}   \right)
	\, ds \right\|  }_{L^2} \, .
\end{align*}
Using Theorem \ref{theoop} and the dispersive estimate \eqref{dispp}, we can bound the above quantity by
\begin{align*}
& \int_1^t s
	{\left\|  
	T_{ \mu_0 (\eta,\s) } \left( e^{is\Lambda} h, e^{is\Lambda} h \right) \right\|}_{L^4}
	{\left\|  e^{it\Lambda} \frac{1}{\grad}  h \right\|}_{L^4} \, ds
\\
& \lesssim \int_1^t s
	{\left\| e^{is\Lambda} h \right\|}_{L^8}^2
	\frac{1}{\sqrt{s}} {\left\| h \right\|}_{L^\frac{4}{3}} \, ds	
\\
& \lesssim \int_1^t s \frac{1}{s^\frac{3}{2}} {\| u \|}_{X}^2
	\frac{1}{\sqrt{s}} {\left\|  \langle x \rangle h \right\|}_{L^2} \, ds
\lesssim  \log t \,  {\| u \|}_{X}^3  \,  .
\end{align*}

By Plancharel, Theorem \ref{theoop},  and \eqref{fractional}  we have
\begin{align*}
{ \left\| \eqref{xHhat2} \right\| }_{L^2} =&
	{\left\|  \int_1^t  \frac{1}{\grad}  e^{-is\Lambda} T_{ \mu_0 (\eta,\s) }
	\left( e^{is\Lambda} h, e^{is\Lambda} h \right) \,  e^{is\Lambda} x h  \, ds \right\| }_{L^2} 
	\\
\lesssim &  \int_1^t  {\left\| T_{ \mu_0 (\eta,\s) } \left( e^{is\Lambda} h, e^{is\Lambda} h \right)
			\right\|}_{L^{24} \cap L\frac{8}{3}} 
			{\left\|  x h \right\|}_{L^2}\, ds	
\\
\lesssim & {\| u \|}_{X} \int_1^t {\left\|  e^{is\Lambda} h \right\|}_{L^\frac{16}{3}}^2 \, s^\g \, ds
	\lesssim {\| u \|}_{X}^3  \int_1^t \frac{1}{s^\frac{5}{4}}  s^\g \, ds  \lesssim {\| u \|}_{X}^3  \,  .
\end{align*}

%

\subsubsection{Estimate of $\sup_t {\| \grad x H (h,h,h) \|}_{H^1}$}
We apply $|\xi| \nabla_\xi$ to $\hat{H}(h,h,h)$, use \eqref{xidxiphi3},
and integrate by parts in time and frequency.
The main contributions after these manipulations  are 
\begin{subequations}
\begin{align}
\label{xidxiH1}
& \iint t \, e^{it\phi(\xi,\eta,\s)}  \frac{\mu_0(\eta,\s)}{|\eta|}
		\what{h}(\s) \what{h}(\eta-\s) \what{Rh}(\xi-\eta) \, d\eta d\s 
\\
\label{xidxiH2}
& \int_1^t \! \iint e^{it\phi(\xi,\eta,\s)}  \frac{ \mu_0(\eta,\s)}{|\eta|}
		\what{h}(\s) \what{h}(\eta-\s) |\xi-\eta|  \nabla_\xi \what{h}(\xi-\eta) \, d\eta d\s \, ds.
\end{align}
\end{subequations}
Using Sobolev's embedding and Theorem \ref{theoop}, 
the first term can be handled as follows:
\begin{align*}
{ \left\| \langle \xi \rangle \eqref{xidxiH1} \right\| }_{L^2} & = t
	{\left\| \frac{1}{\grad} T_{ \mu_0 (\eta,\s) }
	\left( e^{it\Lambda} h, e^{it\Lambda} h \right) \,  e^{it\Lambda} R h  \right\| }_{H^1}
\\
& \lesssim t
	{\left\| \frac{1}{\grad} T_{ \mu_0 (\eta,\s) } \left( e^{it\Lambda} h, e^{it\Lambda} h \right) \right\|}_{W^{1,6}}
	{\left\|  e^{it\Lambda} h \right\|}_{W^{1,3}}
\\
& \lesssim t
	{\left\| T_{ \mu_0 (\eta,\s) } \left( e^{it\Lambda} h, e^{it\Lambda} h \right) \right\|}_{H^1}
	\frac{1}{t^{\frac{1}{3}}} {\| u \|}_{X}
\lesssim t
	{\left\|  e^{it\Lambda} h \right\|}_{W^{1,4}}^2
	\frac{1}{t^{\frac{1}{3}}} {\| u \|}_{X}
\lesssim  \frac{1}{t^\frac{1}{3}}  {\| u \|}_{X}^3  \,  .
\end{align*}
The second term is similar to \eqref{xHhat2} and can be estimated in the same way,
since the presence of an extra derivative does not cause any harm.

\subsubsection{Estimate of $\sup_t t^{-a} {\left\| \grad {|x|}^2 H (h,h,h) \right\|}_{L^2}$
and $\sup_t t^{-b} {\left\| \grad^2 {|x|}^2 H (h,h,h) \right\|}_{L^2}$}
We need to apply $\nabla_\xi$ to \eqref{xidxiH1} and \eqref{xidxiH2} and
estimate the $L^2( {\langle \xi \rangle}^2 d\xi)$-norm of the resulting contributions.
Applying $\nabla_\xi$ to \eqref{xidxiH1}  gives as main terms
\begin{subequations}
\begin{align}
\label{xidxi^2H1}
& \iint t^2 \, e^{it\phi(\xi,\eta,\s)}  \frac{\nabla_\xi \phi(\xi,\eta,\s) \mu_0(\eta,\s)}{|\eta|}
		\what{h}(\s) \what{h}(\eta-\s) \what{R h}(\xi-\eta) \, d\eta d\s 
\\
\label{xidxi^2H2}
& \iint t \, e^{it\phi(\xi,\eta,\s)}  \frac{\mu_0(\eta,\s)}{|\eta|}
		\what{h}(\s) \what{h}(\eta-\s) \nabla_\xi \what{h}(\xi-\eta) \, d\eta d\s \, .
\end{align}
\end{subequations}
Applying it to \eqref{xidxiH2}, and using \eqref{dxiphi/eta}, gives:
\begin{subequations}
\begin{align}
\label{xidxi^2H3}
& \int_1^t \! \iint e^{is\phi(\xi,\eta,\s)}  \frac{ \mu_0(\eta,\s)}{|\eta|}
		\what{h}(\s) \what{h}(\eta-\s) |\xi-\eta| \nabla_\xi^2 \what{h}(\xi-\eta) \, d\eta d\s ds
\\
\label{xidxi^2H4}
& \int_1^t \! \iint s \, e^{is\phi(\xi,\eta,\s)}  \mu_0(\xi,\eta) \mu_0(\eta,\s)
		\what{h}(\s) \what{h}(\eta-\s) \nabla_\xi \what{h}(\xi-\eta) \, d\eta d\s ds \, .
\end{align}
\end{subequations}

To estimate \eqref{xidxi^2H1} we use \eqref{dxiphi/eta}, so that
\begin{subequations}
\begin{align}
\nn
{ \left\| \langle \xi \rangle \eqref{xidxi^2H1} \right\| }_{L^2} & = t^2
	{\left\|  T_{\mu_0 (\xi,\eta)} \left( T_{ \mu_0 (\eta,\s) }
	\left( e^{it\Lambda} h, e^{it\Lambda} h \right) , \frac{1}{\Lambda} e^{it\Lambda} h  \right)
	\right\| }_{H^1}
\\
&
\label{xidxi^2H11}
\lesssim t^2 {\left\|  T_{ \mu_0 (\eta,\s) }
	\left( e^{it\Lambda} h, e^{it\Lambda} h \right) \right\| }_{L^p}
	{ \left\|  \frac{1}{\Lambda} e^{it\Lambda} h  \right\| }_{W^{1,q}}
\\
\label{xidxi^2H12}
& + t^2 {\left\|  T_{ \mu_0 (\eta,\s) }
	\left( e^{it\Lambda} h, e^{it\Lambda} h \right) \right\| }_{W^{1,r}}
	{ \left\|  \frac{1}{\Lambda} e^{it\Lambda} h  \right\| }_{L^s}
\end{align}
\end{subequations}
where $\frac{1}{p} + \frac{1}{q} = \frac{1}{r} + \frac{1}{s} = \frac{1}{2}$
and $2 < p,q,r,s < \infty$.
Since \eqref{xidxi^2H11} can be bounded easily by choosing $p = q = 4$,
we just show details for the estimate of \eqref{xidxi^2H12}.
Choosing $r = \frac{24}{7}$, using theorem \ref{theoop}, and the dispersive estimate,
we can bound:
\begin{align*}
\eqref{xidxi^2H12} \lesssim & \, t^2 {\left\|  T_{ \mu_0 (\eta,\s) }
	\left( e^{it\Lambda} h, e^{it\Lambda} h \right) \right\| }_{W^{1,\frac{24}{7}}}
	{ \left\|  \frac{1}{\Lambda} e^{it\Lambda} h  \right\| }_{L^\frac{24}{5}}
\\
\lesssim & \, t^2
	{\left\|e^{it\Lambda} h \right\|}_{W^{1,4}} {\left\|e^{it\Lambda} h \right\|}_{L^{24}}
	\, \frac{1}{ t^\frac{7}{12} } {\left\|  {\Lambda}^\frac{1}{6}  h \right\|}_{L^\frac{24}{19}}
\\
\lesssim & \, t^2
	\frac{1}{\sqrt{t}}  {\left\| \grad x h \right\|}_{H^1} \, 
	\frac{1}{t^\frac{11}{12}} {\| u \|}_{X} \,
	\frac{1}{ t^\frac{7}{12} } {\left\| x h \right\|}_{H^1}
\lesssim 
	t^\g {\| u \|}_{X}^3
\lesssim
	 t^b {\| u \|}_{X}^3  \,   .
\end{align*}

The term \eqref{xidxi^2H2} is similar to \eqref{xHhat2} except
that it has an extra factor of $t$, but no time integration.
Since these two facts compensate each other, we can estimate \eqref{xidxi^2H2}
in a similar fashion,
obtaining ${ \left\| \langle \xi \rangle \eqref{xidxi^2H2} \right\| }_{L^2} \lesssim {\| u \|}_{X}^3$.
We skip the details.

To estimate \eqref{xidxi^2H3} we note that by 
 Plancharel, Theorem \ref{theoop},  and since $a \leq \frac{1}{3}$ we have
\begin{align*}
{ \left\| \langle \xi \rangle \eqref{xidxi^2H3} \right\| }_{L^2}  =  &{\left\|
	\int_1^t    e^{-is\Lambda} \frac{1}{\grad}  T_{ \mu_0 (\eta,\s) } 
	\left( e^{is\Lambda} h, e^{is\Lambda} h \right) \, e^{is\Lambda} \grad {|x|}^2 h   \, ds
	\right\|  }_{{H}^1}  
\\
\lesssim 
& \int_1^t {\left\|  \frac{1}{\grad} T_{ \mu_0 (\eta,\s) }
	\left( e^{is\Lambda} h, e^{is\Lambda} h \right) \right\| }_{ W^{1,\frac{1}{\d}} }
	{\left\| \grad {|x|}^2 h \right\| }_{H^1} ds
\\
\lesssim 
&  \int_1^t  {\left\|  T_{ \mu_0 (\eta,\s) }
	\left( e^{is\Lambda} h, e^{is\Lambda} h \right) \right\| }_{L^{ \frac{3}{1 + 3\d} } 
									\cap L^\frac{1}{\d}}
	{\left\| \grad {|x|}^2 h \right\| }_{H^1} ds
\\
\lesssim
&  \int_1^t   {\left\| e^{is\Lambda} h \right\| }^2_{ L^\frac{6}{1 + 3\d} }    \, s^a {\| u \|}_{X} ds
\lesssim
  \int_1^t  \frac{1}{s^\frac{4}{3}} \, s^{2\d} \, s^ads  {\| u \|}_{X}^3   \lesssim   t^b {\| u \|}_{X}^3 ,
\end{align*}

Finally, we need to estimate the contribution coming from \eqref{xidxi^2H4}.
By Plancharel's Theorem we see that
\begin{align*}
{ \left\| \eqref{xidxi^2H4} \right\| }_{L^2}  \lesssim  &
	\int_1^t  s {\left\|    T_{ \mu_0 (\xi,\eta) } \left(  T_{ \mu_0 (\eta,\s) } 
	\left( e^{is\Lambda} h, e^{is\Lambda} h \right) , e^{is\Lambda} |x| h  \right) \right\|  }_{L^2}  \, ds
\\\lesssim 
&  \int_1^t  s{\left\|  T_{ \mu_0 (\eta,\s) }
	\left( e^{is\Lambda} h, e^{is\Lambda} h \right) \right\| }_{ L^4 }
	{\left\| e^{is\Lambda} {|x|} h \right\| }_{L^4} ds
\\
\lesssim
& \int_1^t  s {\left\|  e^{is\Lambda} h \right\| }^2_{L^8}
	\frac{1}{\sqrt{s}} {\left\| \grad {\langle x \rangle}^2 h \right\| }_{H^1} ds
\lesssim 
\int_1^t   \frac{1}{s} \, s^{a} ds   {\| u \|}_{X}^3  \lesssim  t^a {\| u \|}_{X}^3  \, .
\end{align*}
For the $L^2( |\xi|^2 d\xi)$-norm we have
\begin{align*}
{ \left\| | \xi | \eqref{xidxi^2H4} \right\| }_{L^2} & \lesssim  
	\int_1^t  s {\left\|    T_{ \mu_0 (\xi,\eta) } \left(  T_{ \mu_0 (\eta,\s) } 
	\left( e^{is\Lambda} h, e^{is\Lambda} h \right) , e^{is\Lambda} |x| h  \right) \right\|  }_{\dot{H}^1}   ds	\, .
\end{align*}
If the derivative in the $\dot{H}^1$ norm falls on the term $e^{is\Lambda} |x| h$,
we can proceed by performing an $L^8 \times L^8 \times L^4$ estimate.
This will yield the bound ${ \left\| | \xi | \eqref{xidxi^2H4} \right\| }_{L^2} \lesssim t^b$.
If instead the derivative hits $T_{ \mu_0 (\eta,\s) } \left( e^{is\Lambda} h, e^{is\Lambda} h \right)$,
we estimate as follows:
\begin{align*}
\int_1^t  s  {\left\|  \Lambda T_{ \mu_0 (\eta,\s) }
	\left( e^{is\Lambda} h, e^{is\Lambda} h \right) \right\| }_{ W^{1,\frac{1}{\d}} }
	{\left\| {|x|} h \right\| }_{L^2} ds
\lesssim
&  \int_1^t  s  {\left\|  e^{is\Lambda} h \right\| }_{ W^{2,\frac{2}{\d}} }
	{\left\|  e^{is\Lambda} h \right\| }_{ L^{\frac{2}{\d}} }
	{\left\| |x| h \right\| }_{L^2}ds
\\
\lesssim 
& \int_1^t  s   \frac{1}{s^{1-\d}} s^{\d_{N}} {\| u \|}_{X} \frac{1}{s^{1-\d}}  {\| u \|}_{X} s^\g  {\| u \|}_{X}ds
\\
\lesssim
&\int_1^t    \frac{1}{s} \, s^{b} \, ds {\| u \|}_{X}^3    \lesssim  t^b  {\| u \|}_{X}^3 \, ,
\end{align*}
provided $2\d + \d_{N} + \g \leq b$. This  concludes the proof of weighted estimates on $H(h,h,h)$.

\subsection{$L^\infty$ estimate on $H(h,h,h)$}
The $L^\infty$ decay for $e^{it\grad} H(h,h,h)$ can be proven by interpolation of weighted $L^2$-norms as we explain below.
By \eqref{disp0} we just need to bound ${\|  {|x|}^{\frac{3}{2}+} \grad^2 H \|}_{L^2}$.
By commuting the weight $x$ and the derivative $\grad$, is is enough to bound
${\|  {|x|}^{\frac{1}{2}+} H \|}_{L^2} + {\|  {|x|}^{\frac{1}{2}+} \grad |x| \grad H \|}_{L^2}$.
We just focus on the second quantity as the first one can be bounded in an easier fashion.
In the previous paragraph we saw that
$$
x \grad  H (h,h,h) \sim \F^{-1} \eqref{xidxiH1}  +  \F^{-1} \eqref{xidxiH2}
$$
We then estimated 
$$
\left\{
\begin{array}{l}
{\| \langle \xi \rangle \eqref{xidxiH1} \|}_{L^2} \lesssim \frac{1}{t^\frac{1}{3}} {\| u \|}_X^3
\\ 
{\| \langle \xi \rangle \nabla_\xi \eqref{xidxiH1} \|}_{L^2} 
\sim {\| \langle \xi \rangle  \eqref{xidxi^2H1} \|}_{L^2} \lesssim \eqref{xidxi^2H12}
\lesssim t^b {\| u \|}_X^3 \, .
\end{array}
\right.
$$
Interpolating between the above two bounds we obtain
\begin{equation}
\label{LinftyH1}
{\| {|x|}^{\frac{1}{2}+} \grad \F^{-1} \eqref{xidxiH1} \|}_{L^2} \lesssim {\| u \|}_X^3 \, . 
\end{equation}

To show the analogous bound for $\F^{-1} \eqref{xidxiH2}$ we can use 
an interpolation argument similar to the one adopted in section \ref{Linftyh}.
Namely, we first write
$$
\F^{-1} \eqref{xidxiH2} = \int_1^t A(s) \, ds \, ,
$$
with the natural definition for $A(s)$.
Then it is easy to see that 
$$
{\| A \|}_{H^1} \lesssim \frac{1}{s^{\frac{4}{3}-}} {\| u \|}_X^3 \, .
$$
Furthermore, in the course of the estimates performed in the previous paragraph
on $|\xi| \nabla_\xi \eqref{xidxiH2} \sim |\xi| \eqref{xidxi^2H3}  +  |\xi| \eqref{xidxi^2H4}$,
we already showed
$$
{\| x \grad A \|}_{L^2} \lesssim \frac{1}{s} s^b {\| u \|}_X^3 \, .
$$
By interpolating these last two bounds we get
$$
{\| {|x|}^{\frac{1}{2}+} \grad \F^{-1} \eqref{xidxiH2} \|}_{L^2} 
\lesssim  \int_1^t {\| {|x|}^{\frac{1}{2}+} \grad A \|}_{L^2}  \, ds
\lesssim  \int_1^t {\| \grad A \|}_{L^2}^{\frac{1}{2}-}  {\| x \grad A \|}_{L^2}^{\frac{1}{2}+}  \, ds
\lesssim  {\| u \|}_X^3 \, ,
$$
which together with \eqref{LinftyH1} implies the desired $L^\infty$ bound on $H(h,h,h)$.

\subsection{Estimates on $H(g,f,f)$}
The term $H(g,f,f)$ is given by the following quartic expression:
\begin{equation}
\label{Hgff} 
\int_1^t e^{-is\grad} \grad^{-1} T_{\mu_0 (\eta,\s)}  \left( e^{\pm is \grad} f \, e^{\pm is \grad}  \grad^{-1} f
		\, , e^{\pm is \grad} f \right) \, e^{is \grad} f \, ds \, .
\end{equation}
We skip the estimate of ${\| x \eqref{Hgff} \|}_{H^1}$ as this is easy to obtain,
and directly move on to estimate ${\| \grad {|x|}^2 \eqref{Hgff} \|}_{H^1}$.
Applying $\nabla_\xi^2$ to $\hat{H}(g,f,f)$, and using \eqref{dxiphi/eta},
produces two main contributions, which expressed in real space are of the form:
\begin{subequations}
\begin{align}
\label{Hgff1}
& \int_1^t s^2 e^{-is\grad}  T_{\mu_0 (\xi,\eta)} \left[
	T_{\mu_0 (\eta,\s)}  \left( 
		u \, \grad^{-1} u \, , u \right) , \grad^{-1} u \right] \, ds \,
\\ 
\label{Hgff2}
& \int_1^t e^{-is\grad} \grad^{-1} T_{\mu_0 (\eta,\s)}  \left( u \,  \grad^{-1} u
		\, , u \right) \, e^{is \grad} {|x|}^2 f \, ds \, .
\end{align}
\end{subequations}
We now estimate the $\dot{H}^1$-norm of the above terms.
Using Theorem \ref{theoop} as usual, and \eqref{u/gradLp}, we can estimate
\begin{align*}
{\| \eqref{Hgff1} \|}_{\dot{H}^1} & \lesssim
\int_1^t  s^2  {\left\|  T_{\mu_0 (\eta,\s)}  \left( 
		u \, \grad^{-1} u \, , u \right) \right\| }_{ W^{1,4} }
	{\left\| \grad^{-1} u \right\| }_{W^{1,4}} \, ds
\\
& \lesssim  \int_1^t  s^2  {\left\|  u \right\| }^2_{ W^{2,\frac{1}{2\d}} }
	{\left\| \grad^{-1} u \right\| }_{W^{1,4}}^2 \, ds
\\
& \lesssim 
\int_1^t  s^2   \frac{1}{s^{2-2\d}} s^{2\d_{N}} {\| u \|}^2_{X} \frac{1}{s} s^{2\g} {\| u \|}^2_{X}  \, ds
\lesssim  t^{2\d + 2\d_N + 2\g}  {\| u \|}_{X}^4 \, ,
\end{align*}
To estimate \eqref{Hgff2} we distinguish two different cases:
one where the derivative in the $\dot{H}^1$-norm falls on $e^{is \grad} {|x|}^2 f$
and another one when it falls on 
$\grad^{-1} T_{\mu_0 (\eta,\s)}  \left( u \,  \grad^{-1} u \, , u \right)$.
In the first case we can estimate:
\begin{align*}
& {\left\| \int_1^t e^{-is\grad} \grad^{-1} T_{\mu_0 (\eta,\s)}  \left( u \,  \grad^{-1} u
		\, , u \right) \, e^{is \grad} \grad {|x|}^2 f \, ds \right\|}_{L^2}
\\
& \lesssim \int_1^t  {\left\|  \grad^{-1} T_{\mu_0 (\eta,\s)}  \left( 
		u \, \grad^{-1} u \, , u \right) \right\| }_{ W^{1,\frac{1}{\d}} }
		{\left\| \grad {|x|}^2 f \right\| }_{L^2} \, ds
\\
& \lesssim  \int_1^t   {\left\|  T_{\mu_0 (\eta,\s)}  \left( 
		u \, \grad^{-1} u \, , u \right)  \right\| }_{ W^{1,\frac{3}{1 + 3\d}} }
	{\left\| \grad {|x|}^2 f \right\| }_{L^2} \, ds
\\
& \lesssim 
\int_1^t   {\left\|  u \, \grad^{-1} u \right\| }_{ W^{1,\frac{6}{2 + 3\d}} }
	{\left\|  u \right\| }_{ W^{1,\frac{2}{\d}} }
	{\left\| \grad {|x|}^2 f \right\| }_{L^2} \, ds
\lesssim 
\int_1^t   {\left\|  u \right\| }^2_{ W^{1,\frac{2}{\d}} }  {\| \grad^{-1} u \|}_{W^{1,3}} 
	{\left\| \grad {|x|}^2 f \right\| }_{L^2} \, ds
\\
&
\lesssim 
\int_1^t  \frac{1}{s^{2-2\d}} s^{2\d_N} {\| u \|}_{X}^2 s^\g {\| u \|}_{X}  s {\| u \|}_{X} \, ds
\lesssim t^{2\d + 2\d_N + \g}  {\| u \|}_{X}^4 \, .
\end{align*}
In the second case we use again Theorem \ref{theoop} and \eqref{u/gradLp} to obtain
\begin{align*}
& {\left\| \int_1^t e^{-is\grad} T_{\mu_0 (\eta,\s)}  \left( u \,  \grad^{-1} u
		\, , u \right) \, e^{is \grad} {|x|}^2 f \, ds \right\|}_{L^2} & 
\\
& \lesssim \int_1^t  {\left\| T_{\mu_0 (\eta,\s)}  \left( 
		u \, \grad^{-1} u \, , u \right) \right\| }_{L^3}
		{\left\| {|x|}^2 f \right\| }_{L^6} \, ds
\lesssim 
\int_1^t   {\left\|  u \right\| }^2_{ W^{1,\frac{2}{\d}} }  {\| \grad^{-1} u \|}_{W^{1,3}} 
	{\left\| \grad {|x|}^2 f \right\| }_{L^2} \, ds
\\
& \lesssim 
\int_1^t  \frac{1}{s^{2-2\d}} s^{2\d_N} {\| u \|}_{X}^2 s^\g {\| u \|}_{X}  s {\| u \|}_{X} \, ds
\lesssim t^{2\d + 2\d_N + \g}  {\| u \|}_{X}^4 \, .
\end{align*}
For appropriate choices of $\d,\g$ and $N$,
and since adding a derivative only costs a $t^{\d_N}$ factor, the above bounds guarantee
${\| \grad \eqref{Hgff1} \|}_{\dot{H}^1}, {\| \grad \eqref{Hgff2} \|}_{\dot{H}^1} \lesssim t^b  {\| u \|}_{X}^4$.
This give that ${\| \grad {|x|}^2 \eqref{Hgff} \|}_{H^1} \lesssim t^b  {\| u \|}_{X}^4$,
which is an even stronger bound that what is needed.

Finally, the $L^\infty$ estimate on $e^{it\grad} H(g,f,f)$ can be obtained by interpolating 
weighted $L^2$-norms as we have already done in previous sections.

\appendix
\section{Tools from harmonic analysis\label{appendix}}
In this appendix we gathered some elementary results and some multilinear estimates from Harmonic analysis,
as well as the proof of Theorem \ref{theoop}.
We start by recalling some dispersive properties of the linear wave propagator $e^{it\Lambda}$:
\begin{lem}
In three space dimensions the following dispersive estimates hold:
\begin{align}
\label{disp0}
{\left\| e^{it\Lambda} f \right\|}_{L^\infty} & \lesssim \frac{1}{t} 
	{\| {\langle x \rangle}^{\frac{3}{2}+} \grad^2 f \|}_{L^2} 
\\
\label{disp}
{\left\| e^{it\Lambda} f \right\|}_{L^\infty} & \lesssim \frac{1}{t}
	{\left\| | x | \grad^2 f \right\|}^\frac{1}{2}_{L^2} {\left\| {| x |}^2 \grad^2 f \right\|}^\frac{1}{2}_{L^2}  \, .
\end{align}
Furthermore, for $2 \leq p < \infty$
\begin{equation}
\label{dispp}
{\left\| e^{it\Lambda} f \right\|}_{L^p} \lesssim  
                  \frac{1}{t^{1-\frac{2}{p}}} {\|\Lambda^{2-\frac{4}{p}} f \|}_{L^{p^\p}} \, .
\end{equation}
where $p^\p$ is the dual H\"{o}lder exponent of $p$.
\end{lem}

\proof
The first two inequalities are a consequence of the standard $L^1 \rightarrow L^\infty$ estimate
for the linear wave equation applied to the propagator $e^{it\Lambda}$:
\begin{equation*}
{\left\| e^{it\Lambda} f \right\|}_{L^\infty} \lesssim \frac{1}{t} 
		\left[ {\| f \|}_{\dot{W}^{2,1}} + \ {\| \grad f \|}_{\dot{W}^{1,1}} \right] \, ,
\end{equation*}
see for instance \cite{SS}, together with the fact that
${\left\| g \right\|}_{L^1} \leq {\left\| | x | g \right\|}^\frac{1}{2}_{L^2} {\left\| {| x |}^2 g \right\|}^\frac{1}{2}_{L^2}$.
The estimate \eqref{dispp} is a consequence of the dispersive estimate
\begin{equation*}
{\left\| e^{it\Lambda} f \right\|}_{\dot{B}^0_{p,2}} \lesssim \frac{1}{t^{1-\frac{2}{p}}} 
		{\| f \|}_{\dot{B}^{2-\frac{4}{p}}_{p^\p,2}} \, ,
\end{equation*}
and basic relations between Besov-norms and $L^p$-norms  see again \cite{SS} $\, _\Box$

We also recall two standard inequalities about fractional integration:
\begin{align}
\label{fractional}
{\left\|  \frac{1}{{\Lambda}^\alpha} \, f  \right\|}_{L^q} & \lesssim  {\| f \|}_{L^p} 
\quad \mbox{for \, $1 < p,q < \infty$ \, and \, $\alpha = \frac{3}{p} - \frac{3}{q}$} \, ;
\\
\label{fractional2}
{\left\|  \frac{1}{{\Lambda}^\alpha} \, e^{it\Lambda}  f  \right\|}_{L^q}  & \lesssim  {\| f \|}_{L^p} 
\quad \mbox{for \, $1 < p \leq 2 \leq q < \infty$ ,  $\alpha = \frac{3}{p} - \frac{3}{q}$ 
\, and \, $0 < \alpha < \frac{3}{p}$
} \, .
\end{align}



\begin{theo}
\label{theoop}

Let $p,q,r$ be given such that $\frac{1}{r} = \frac{1}{p} + \frac{1}{q}$
and $1 < p,q,r < \infty$. The following hold
\begin{enumerate}
\item[(i)] 
\label{theoop1}
If $m$ belongs to the class $\B_0$
\begin{equation*}
{\| T_{m} (f,g) \|}_{L^r} \lesssim {\| f \|}_{L^p} {\| g \|}_{L^q} \, .
\end{equation*}

\item[(ii)]
\label{theoop2}
If $m$ belongs to the class $\B_s$ for $s \geq 0$ and $k$ is an integer, then
\begin{equation*}
{\left\| \grad^k T_{m} (f,g) \right\|}_{L^r} \lesssim {\| f \|}_{W^{s+k,p}} {\| g \|}_{L^{q}} 
				+ {\| f \|}_{L^{p}} {\| g \|}_{W^{s+k,q}} \, .
\end{equation*}

\item[(iii)]
\label{theoop3}
If $m$ belongs to the class $\B_s$ and $M > 3$, then
\begin{equation*}
{\left\| \grad^k T_{m} (f,g) \right\|}_{L^2} \lesssim {\| f \|}_{H^{s+k}} {\| g \|}_{W^{1,M}}
				+ {\| f \|}_{W^{1,M}} {\| g \|}_{H^{s+k}}  \, .
\end{equation*}
\end{enumerate}
\end{theo}

\proof

\label{prooftheoop}
Point (i) is, essentially, equivalent to Theorem C.1 in \cite[Appendix C]{GMS2}.
(ii) is a consequence of the proof of (i).
(iii) can also be easily obtained from the proof of (i) and Sobolev's embedding;
it is a substitute for the lack of an $L^2 \times L^\infty$ estimate.
For completeness we provide the proof of (i) and (ii) below.

Let $m$ be a symbol in $\B_0$. Away from the coordinate axes $\{ \xi = 0 \} \cup \{ \eta = 0 \} \cup \{ \xi - \eta= 0 \}$
the Coifman-Meyer theorem \cite{CM} applies to give the desired boundedness result.
Let us now consider the case $|\eta| \ll |\xi|,|\xi-\eta| \sim 1$.
Notice that this is the only case we really need to treat since the other 
cases can be reduced to this one by duality.
Therefore, we can assume
\begin{equation*}
T_m (f,g) = \sum_j T_m (P_{<j-100}f, P_{<j} g) \, .
\end{equation*}
From the definition of the class $\B_0$, close to the $\eta = 0$ axes we can 
assume $m = \A \left( |\eta|, \frac{\eta}{|\eta|}, \xi \right)$, for some smooth function $\A$.
By homogeneity $m = \A \left( \frac{|\eta|}{|\xi|}, \frac{\eta}{|\eta|}, \frac{\xi}{|\xi|} \right)$.
Expanding this expression in $\frac{|\eta|}{|\xi|}$ gives
\begin{equation*}
m(\xi,\eta) = \sum_{k=1}^L \frac{ {|\eta|}^k }{ {|\xi|}^k } m_k \left(\frac{\eta}{|\eta|}, \frac{\xi}{|\xi|} \right)
\quad + \quad \mbox{remainder}
\end{equation*}
By assumption on $\A$, the symbols $m_k$ are smooth.
Moreover, if we take $L$ large enough, the singularity of the remainder at $\eta = 0$ 
becomes so weak that the remainder satisfies estimates of Coifman-Meyer type. 
We can then disregard it in what follows.

Expanding $m_k$ in spherical harmonics (denoted by $Z_l$, $l \in \mathbb{N}$) yields
\begin{equation*}
m(\xi,\eta) = \sum_{k=1}^L \sum_{l,l^\p} a_{k,l,l^\p}
\frac{ {|\eta|}^k }{ {|\xi|}^k } Z_l \left(\frac{\eta}{|\eta|} \right) Z_{l^\p} \left( \frac{\xi}{|\xi|} \right) \, .
\end{equation*}
By the Mihlin-H\"{o}rmander multiplier theorem, the operators associated
to the symbols $Z_l \left( \frac{\eta }{| \eta |} \right)$ and $Z_l \left( \frac{\xi}{| \xi |} \right)$
are bounded on Lebesgue spaces with bounds growing polynomially in $l$;
on the other hand, since $m$ is smooth, the coefficients $a_{k,l,l^\p}$ 
decay faster than any polynomial in $(l,l^\p)$. 
We can then ignore the summation over $(l,l^\p)$ and the finite summation over $k$.
Thus, matters reduce to bound
\begin{equation*}
\sum_j \grad^{-k} \left(  P_{< j-100} \grad^k f P_j g  \right) \, .
\end{equation*}
Using the Littlewood-Paley square and maximal function estimates we finally obtain
\begin{align*}
& {\left\| \sum_j \grad^{-k} \left(  P_{< j-100} \grad^k f P_j g  \right) \right\|}_{L^r} 
\lesssim 
{\left\| {\left[ \sum_j 2^{-2jk} {\left(  P_{< j-100} \grad^k f P_j g  \right)}^2  \right]}^\frac{1}{2} \right\|}_{L^r} 
\\
& \lesssim {\left\|  \sup_j \left| 2^{-jk} P_{< j-100} \grad^k f \right| \right\|}_{L^p}
{\left\| {\left[ \sum_j {\left( P_j g  \right)}^2  \right]}^\frac{1}{2} \right\|}_{L^q}  
\lesssim {\| f \|}_p  {\| g \|}_q \, \,  _\Box
\end{align*}



\addcontentsline{toc}{section}{Bibliography}
\bibliographystyle{plain}

\def\cprime{$'$}

\end{document}